\documentclass[12pt,onecolumn,draftcls]{IEEEtran}

\usepackage{amsmath,amsthm,amssymb,float,graphicx,epstopdf,geometry,bm}
\usepackage{bbm,bbding,amssymb,pifont,mathrsfs,amsfonts,graphicx,mathtools,color}%
\usepackage{amscd,array,enumerate,dsfont,texdraw,tikz,multicol,bbm}

\newtheorem{lemma}{Lemma}
\newtheorem{assumption}{Assumption}
\newtheorem{proposition}{Proposition}
\newtheorem{remark}{Remark}
\newtheorem{corollary}{Corollary}
\newtheorem{theorem}{Theorem}

\ifCLASSINFOpdf
\else
\fi

\begin{document}
	
\title{Convergence Properties of the Distributed  Projected Subgradient Algorithm over General Graphs}

\author{Weijian~Li$^*$, Zihan~Chen$^\dag$, Youcheng~Lou$^\ddag$, and Yiguang~Hong$^\S$

\thanks{$^*$W.~Li is with the Department of Automation, University of Science and Technology of China, Hefei, 230027, Anhui, China,
e-mail: \texttt{ustcwjli@mail.ustc.edu.cn}.
			
$^\dag$Z.~Chen is with Huawei Technologies Co., Ltd., Beijing, 100095, China e-mail: \texttt{chenzihan9@huawei.com}.
			
$^\ddag$Y. Lou is with MDIS, Academy of Mathematics and Systems Science, Chinese Academy of Sciences, Beijing, 100190, China, 
e-mail: \texttt{louyoucheng@amss.ac.cn}.
			
$^\S$Y.~Hong is with the
Y.~Hong is with the
Department of Control Science and Engineering, Tongji University, Shanghai, 201804, China,
and also with Shanghai Research Institute for Intelligent Systems, Shanghai, 201203, China,
e-mail: \texttt{yghong@iss.ac.cn}.
}}
	
\maketitle

\begin{abstract}
In this paper, we revisit a well-known distributed projected subgradient algorithm which aims to minimize a sum of cost functions with a common set constraint.
In contrast to most of existing results,  weight matrices of the time-varying multi-agent network are assumed to be more general, i.e., they are only required to be row stochastic instead of doubly stochastic. We focus on analyzing convergence properties of this algorithm under general graphs.
We first show that there generally exists a graph sequence such that the algorithm is not convergent when the network switches freely within finitely many general graphs.
Then  to guarantee the convergence of this algorithm under any uniformly jointly strongly connected general graph sequence, we provide a necessary and sufficient condition,  i.e., the intersection of optimal solution sets to all local optimization problems is not empty.
Furthermore, we surprisingly find that the algorithm is convergent for any periodically switching general graph sequence, and the converged solution minimizes a weighted sum of local cost functions, where the weights depend on the Perron vectors of some product matrices of the underlying periodically switching  graphs. 
Finally, we consider a slightly broader class of quasi-periodically switching graph sequences, and show that the algorithm is convergent for any quasi-periodic  graph sequence if and only if the network switches between only two graphs.
\end{abstract}


\section{Introduction}
In the past decade,
distributed convex optimization has received intensive research attention, motivated by its broad applications in various areas including distributed estimation \cite{kar2012distributed}, resource allocation \cite{molzahn2017survey}, and machine learning \cite{boyd2011distributed}.
The basic idea is that in a multi-agent network, all agents cooperate to solve an optimization problem with their local cost functions, constraints and neighbors' states.
A variety of distributed algorithms have been designed for different formulations.
By performing local averaging operations and taking subgradient descent steps, a distributed subgradient algorithm was proposed for unconstrained distributed optimization problems in \cite{nedic2009distributed}.
Then a projected subgradient algorithm was developed to deal with a common set constraint in \cite{nedic2010constrained}. 
The case of nonidentical set constraints was studied in \cite{lin2016distributed,mai2019distributed}.
Following that, by combining projected subgradient methods and primal-dual ideas, distributed primal-dual subgradient algorithms were designed to minimize a sum of local cost functions with set constraints, local inequality and equality constraints \cite{zhu2011distributed}.
Moreover, distributed primal-dual algorithms have also been explored for coupled constraints \cite{cherukuri2016initialization,zeng2018distributed}.

For a distributed algorithm, the communication topology of the multi-agent network has a great effect on its convergence \cite{nedic2018network}.
Since the pioneering work for distributed optimization in \cite{nedic2009distributed},
weight-balanced graphs have been widely employed to design distributed algorithms \cite{gharesifard2013distributed,
yuan2016convergence,zeng2018distributed,
zhu2011distributed}.
In \cite{zeng2018distributed}, undirected and connected graphs were adopted for distributed primal-dual algorithms.
Based on saddle-point dynamics, a continuous-time algorithm was proposed under  strongly connected and weight-balanced digraphs in \cite{gharesifard2013distributed}.
Time-varying weight-balanced graphs have also been utilized for the distributed design \cite{yuan2016convergence,zhu2011distributed}.

Plenty of results on distributed optimization were obtained under weight-balanced graphs because there usually existed a common Lyapunov function to	facilitate the convergence analysis \cite{lou2014approximate,nedic2009distributed,
olshevsky2008nonexistence}.
Furthermore, most of (sub)gradient-based algorithms could achieve an optimal solution under weight-balanced graphs because Perron vectors of weight matrices were with identical entries \cite{lin2016distributed,mai2019distributed,
nedic2009distributed}.
However, a weight-balanced graph requires the in-degree of each node being equal to its out-degree, and is not always practical in real applications \cite{nedic2018network}. 
For instance, if agents use broadcast-based communications in a wireless network, they neither know their out-neighbors nor are able to adjust their outgoing weights. Thus, the weight-balance condition is difficult to be guaranteed in this case \cite{mai2019distributed}.
To overcome the difficulty, some new mechanisms have been explored to design distributed algorithms under weigh-unbalanced graphs \cite{mai2019distributed,pu2020push,tsianos2012push,tsianos2012distributed}.
In \cite{tsianos2012distributed}, a reweighting technique was proposed for fixed graphs with known Perron vectors.
By combining the dual averaging algorithm with the push-sum consensus protocol, a distributed push-sum method was developed in \cite{tsianos2012push}.
Requiring a row stochastic and a column stochastic matrices, a distributed push-pull algorithm was designed in \cite{pu2020push}, where for individual agent, the gradient was pushed to its neighbors, and the decision variable was pulled from its neighbors.
In \cite{lou2015nash,mai2019distributed}, heterogeneous stepsizes were adopted to balance the graph.

Some researches have been focused on developing distributed algorithms under weight-unbalanced graphs.
However, how do unbalanced networks affect the performance of a distributed algorithm?
Note that it is an important problem because it can provide us with a better understanding of existing works, and moreover, guide us to design effective distributed algorithms.

In this paper, we revisit a well-known distributed projected subgradient algorithm, proposed in \cite{nedic2010constrained},  to minimize a sum of (nonsmooth) cost functions with a common set constraint. Compared with existing results in \cite{liu2017convergence, nedic2010constrained}, we assume the time-varying communication network being general (weight matrices of the network are row stochastic instead of doubly stochastic), and focus on analyzing convergence properties of this algorithm. Our main contributions are summarized as follows.

\begin{itemize}
\item We show that there generally exists a graph sequence such that the algorithm is not convergent if the time-varying network switches freely within finitely many general graphs. 
\item  To guarantee the convergence of this algorithm for any uniformly jointly strongly connected  general graph sequence, we provide a necessary and sufficient condition, namely, the intersection of optimal solution sets to all local optimization problems is not empty.
\item We find that the algorithm is convergent for any periodically switching general graph sequence, and moreover, the converged solution minimizes a weighted sum of the local cost functions.
In addition, we relax the periodic condition slightly, and define a broader class of quasi-periodic graph sequences. We show that the algorithm is always convergent for any quasi-periodic graph sequence if and only if the network switches between two graphs. 
\end{itemize}

The remainder of this paper is organized as follows. Some preliminary knowledge related to convex analysis and graph theory is introduced in Section 2, and then the problem is formulated in Section 3.
Our main results are presented Sections 4, while their rigorous proofs are provided in Section 5.
Following that, illustrative examples are carried out in Section 6. 
Finally, concluding remarks are given in Section 7.

\textbf{Notations:}
Let $\mathbb R$, $\mathbb R^m$
and $\mathbb R^{m\times n}$
be the set of real numbers,
the set of $m$-dimensional real column vectors, and the set of $m$-by-$n$ dimensional real matrices, respectively. 
Let $\mathbb N$ be the set of nonnegative integers.
Vectors are column vectors by default.
$x'$ stands for the transpose of vector $x$.
$[A]_{ij}$ means the $(i,j)$-th entry of matrix $A$.
The Euclidean norm of $x,y\in \mathbb R^m$ is defined by $x'y$.
Let $|\cdot|$, $|\cdot|_1$ be the Euclidean norm and $l_1$-norm of a vector, respectively.
Denote ${\rm dist}(x,\Omega)$ as the distance
from a point $x$ to a set $\Omega$ (that is, ${\rm dist}(x,\Omega)\triangleq \inf_{y\in\Omega} |y-x|$).

\section{Preliminary knowledge}
In this section, we introduce some basic concepts related to convex analysis and graph theory.

\subsection{Convex analysis}
A set $\Omega\subset\mathbb R^m$ is convex if
$\lambda x+(1-\lambda)y \in \Omega$ for all
$x,y\in \Omega$ and $\lambda\in [0,1]$.
For a closed convex set $\Omega \subset \mathbb R^m$, we define $P_{\Omega}(\cdot): \mathbb R^m \rightarrow \Omega$ as a projection operator which maps 
$x\in\mathbb R^m$ to a unique point $P_{\Omega}(x)$ such that
$P_{\Omega}(x)={\rm argmin}_{y\in \Omega}|x-y|$.
Referring to Lemma 1 in \cite{nedic2010constrained},
we have
\begin{equation}\label{pro_ine1}
|P_{\Omega}(x)-y| \le |x-y|,
~\forall x\in\mathbb R^m,
~\forall y\in \Omega,
\end{equation}
and moreover,
\begin{equation}\label{pro_ine2}
|P_{\Omega}(x)-y|^2 \le |x-y|^2 -|x-P_{\Omega}(x)|^2, 
~\forall x\in\mathbb R^m,
~\forall y\in \Omega.
\end{equation}

A function $f: \Omega \rightarrow \mathbb R$
is convex if $\Omega$ is a convex set, and
\begin{equation*}
f(\theta x+(1-\theta)y) \le
\theta f(x)+(1-\theta)f(y),
~\forall x,y \in \Omega,
~\forall \theta \in [0,1].
\end{equation*}
Furthermore, it is strictly convex if the strict inequality holds whenever $x\not= y$ and $\theta\in (0,1)$. 
If $g_f(x)\in\mathbb R^m$ satisfies 
$$f(y)-f(x) \ge (y-x)'g_f(x),$$
then $g_f(x)$ is the subgradient of $f$ at $x$. 
Denoted by $\partial f(x)$ the set of all subgradients
of $f$ at $x$.

\subsection{Graph theory}
The communication topology of a multi-agent network can be modeled by a digraph $\mathcal G(\mathcal V, \mathcal E)$, 
where $\mathcal V$ is the node set, and
$\mathcal E$ is the edge set.
Then a nonnegative weight matrix 
$A=[a_{ij}]$ can be associated with $\mathcal G$,
where $a_{ij}>0$ if and only if 
$(j,i)\in \mathcal E$.
Conversely, a graph $\mathcal G$  can also be associated with a nonnegative  matrix $A$. 
Node $j$ is a neighbor of $i$, and can send information to $i$ if  $a_{ij}>0$.
Suppose that there are self-loops in $\mathcal G$, i.e., $a_{ii}>0$ for all $i \in \mathcal V$.
Graph $\mathcal G$ is said to be weight-balanced if $\sum_{j\in\mathcal V}a_{ij}= \sum_{j\in\mathcal V}a_{ji}$ 
for  $i\in\mathcal V$, 
and is weight-unbalanced otherwise.
A path from $i_1$ to $i_p$ is defined by an edge sequence 
$(i_1,i_2), (i_2,i_3), \dots, (i_{p-1},i_p) \in \mathcal E$
with distinct nodes $i_1,\dots, i_p$.
$\mathcal G$ is strongly connected if there exists at least a path between every pair of nodes.
If a network is time-varying, we denote
$\mathcal G(\mathcal V, \mathcal E(k))$ or $\mathcal G(k)$ as the graph at time slot $k$. Furthermore, the joint
graph over the time interval $[k_1,k_2)$ is given by $\mathcal G([k_1,k_2))\triangleq \mathcal G(\mathcal V, \bigcup_{k\in [k_1,k_2)} 
\mathcal E(k))$.

A vector is said to be stochastic if it is with nonnegative entries and the sum of its entries is $1$. Furthermore, it is also positive if all entries of the vector are positive.
A matrix is row (column) stochastic if all of its row (column) vectors are stochastic, 
and is doubly stochastic if it is both row and column stochastic.
A row stochastic matrix is also sometimes simply called a stochastic matrix.
The following result, 
collected from Lemma $5.3$ in \cite{lou2015nash}, 
addresses the relationship of positive stochastic vectors and stochastic matrices.

\begin{lemma}
\label{sto_mat_lem}
For any positive stochastic vector $\mu \in \mathbb R^n$, there
must be a stochastic matrix 
$A=(a_{ij})\in \mathbb R^{n\times n}$
such that $\mu' A=\mu'$, and moreover,
the graph associated with $A$ is strongly connected.
\end{lemma}

Let $B \in \mathbb R^{n\times n}$ be a stochastic matrix, and $\mathcal G_B$ be the associated graph.
It follows from the Perron-Frobenius theorem
\cite{horn2012matrix} that
there is a unique positive stochastic left eigenvector $\mu(B)$ of $B$ associated with eigenvalue $1$ if $\mathcal G_B$ is strongly connected. We call $\mu(B)$ the Perron vector of $B$.

\section{Formulation and algorithm}
In this section, we formulate the constrained distributed optimization problem, and then revisit a projected subgradient algorithm. Furthermore, we give the problem statement.

Consider a network of $n$ agents connected by a time-varying digraph
$\mathcal G(\mathcal V, \mathcal E(k))$ 
(or simply $\mathcal G(k)$),
where
$\mathcal V=\{1,\dots,n\}$ and 
$\mathcal E(k)\subset\mathcal V \times \mathcal V$.
For each $i\in \mathcal V$, there is a local (nonsmooth) cost function 
$f_i:\mathbb R^m \rightarrow \mathbb R$
and a feasible constraint set 
$X \subset \mathbb R^m$.
All agents cooperate to minimize the global cost function $\sum_{i\in\mathcal V} f_i(x)$ in $X$.
To be strict, the problem can be formulated as
\begin{equation} \label{pro_for}
{\rm min}~~\sum_{i\in \mathcal V} f_i(x),
\quad {\rm s.t.}~~x \in X,
\end{equation}
where $x$ is the decision variable.

Let $x_i(k)$ be the estimation for a solution to (\ref{pro_for}) by agent $i$. Then a distributed algorithm is said to achieve a solution to (\ref{pro_for}) if for any initial condition 
$x_i(0) \in \mathbb R^m$,
$\lim_{k\to \infty} |x_i(k)-x_j(k)|=0$, and moreover, there exists $x^* \in X^*$  such that $\lim_{k\to\infty}|x_i(k)-x^*|=0$, where
$$X^*=\Big\{z~|~z={\rm argmin}_{x\in X}
\sum_{i\in \mathcal V} f_i(x)\Big\}.$$

To ensure the well-posedness of (\ref{pro_for}), we make the following standard assumptions.

\begin{assumption}\label{con_ass}
(Convexity) 
For each $i\in \mathcal V$, $f_i$ is a convex function on an open set containing $X$, 
and  $X$ is a closed convex set.
\end{assumption}

\begin{assumption}\label{lip_ass}
(Boundedness of Subgradients) 
For each $i\in \mathcal V$, the subgradient set of $f_i$ is bounded over $X$, i.e., there exists a scalar $L>0$ such that 
\begin{equation}\label{lip_ine}
|d|\le L, 
~\forall d\in \partial f_i(x),
~\forall x\in X.
\end{equation}
\end{assumption}

\begin{assumption}\label{gra_ass}
(Connectivity)
The graph sequence 
$\mathcal G(k)$ 
is uniformly jointly 
strongly connected (UJSC), i.e., there exist an integer $B>0$
such that the joint graph $\mathcal G([k,k+B))$
is strongly connected for 
$k\in \mathbb N$.
\end{assumption}

\begin{assumption}\label{wei_ass}
(Weight Rule)
\begin{enumerate}[(i)]
\item The weight matrix $A(k)=[a_{ij}(k)]$ associated to $\mathcal G(k)$ is stochastic, i.e.,\\
$\sum_{j\in \mathcal V}a_{ij}(k)=1$ for 
$i\in \mathcal V$ and $k\in \mathbb N$.
\item There is a scalar $0<\eta<1$ such that
$a_{ij}(k) \ge \eta$ if $a_{ij}(k)>0$, and 
$a_{ii}(k) \ge \eta$ for $i,j\in \mathcal V$ and $k\in \mathbb N$.
\end{enumerate}
\end{assumption}

Note that (\ref{pro_for}) is a well-known constrained distributed optimization problem investigated in \cite{liang2019dual,nedic2010constrained,
qiu2016distributed}.
A pioneering distributed algorithm for this problem is the projected subgradient method, which combines an average step with a local projected gradient update step \cite{nedic2010constrained}. The specific form of this algorithm is given by
\begin{equation}\label{alg}
\begin{aligned}
v_i(k)=&
\sum_{j\in\mathcal V}a_{ij}(k)x_j(k),\\
x_i(k+1)=&
P_X(v_i(k)-\alpha_k d_i(k)),
\end{aligned}
\end{equation}
where $d_i(k) \in \partial f_i(v_i(k))$, and $\alpha_k>0$ is the stepsize.
To guarantee the convergence of (\ref{alg}), the following assumption is made \cite{nedic2010constrained}.

\begin{assumption}\label{step_ass}
(Stepsize Rule)
$\sum_{k=0}^\infty\alpha_k=\infty$, 
and moreover,
$\sum_{k=0}^\infty \alpha^2_k<\infty$.	
\end{assumption}

\begin{remark}
In fact, Assumptions \ref{con_ass}-\ref{step_ass} have also been employed in \cite{liu2017convergence,lou2015nash,nedic2010constrained,zhu2011distributed}.
As a comparison, we only suppose that the weight matrix $A(k)$ is row stochastic instead of doubly stochastic, i.e., the communication graph may be weight-unbalanced. Thus, the considered problem is more general.
\end{remark}

The following result, proved in \cite{nedic2010constrained}, shows a convergence property of (\ref{alg}).

\begin{proposition}\label{con:bal}
Let Assumptions \ref{con_ass}-\ref{step_ass} hold. If $A(k)$ is also column stochastic for $k\in \mathbb N$, then algorithm (\ref{alg}) achieves a solution to (\ref{pro_for}).
\end{proposition}

Proposition \ref{con:bal} indicates the convergence of (\ref{alg}) under graphs with  doubly stochastic weight matrices. 
Following that, great efforts have been paid to  develop distributed algorithms over weight-balanced graphs \cite{gharesifard2013distributed, nedic2017achieving, zhu2011distributed}. 
Furthermore, some new mechanisms have  been proposed to replace the weight-balance assumption including the push-sum protocol 
\cite{ tsianos2012push}
and the push-pull method \cite{pu2020push}.
To distinguish with weight-balanced graphs,
we call a graph being general if its weight matrix is only required to be row stochastic.
An interesting question is what the convergence performance of a distributed algorithm is if
the underlying graph is general.
In this paper, taking (\ref{alg}) as a starting point, we explore its convergence under general graphs.
To be specific, we are interested in the following three questions.
\begin{itemize}
\item Is algorithm (\ref{alg}) convergent under general graphs?
\item If it is, which solution does it converge to?  If not, under what condition on the cost functions it is convergent?
\item Is there any class of general graph sequences under which the algorithm (\ref{alg}) is convergent?
\end{itemize}

\section{Main results}
\label{results}
In this section, we present the main results on the convergence  of (\ref{alg}) under general graph sequences.
At the beginning, we show that there generally exists a graph sequence such that (\ref{alg}) is not convergent.
Then we provide a necessary and sufficient condition to guarantee its convergence.
Finally, we establish its convergence under periodic and quasi-periodic graph sequences.

\subsection{Basic results}
Define 
$y(k)=\frac 1n\sum_{i\in\mathcal V} x_i(k)$ as the average of agents' estimations. 
The following lemma, proved in Section \ref{fix:con_lem_pf}, shows consensus results of (\ref{alg}).

\begin{lemma}\label{cons:lem}
Consider algorithm (\ref{alg}). Under Assumptions \ref{lip_ass}-\ref{wei_ass}, the following statements hold.
\begin{enumerate}[(i)]
\item If the stepsize satisfies 
$\lim_{k\to\infty}\alpha_k=0$, then
\begin{equation*}
\lim_{k\to\infty}|x_i(k)-y(k)|=0,
~\forall i\in\mathcal V.
\end{equation*}
\item If the stepsize satisfies
$\sum_{k=0}^\infty \alpha_k^2<\infty$, then
\begin{equation*}
\sum_{k=0}^\infty \alpha_k|x_i(k)-y(k)|<\infty,
~\forall i\in\mathcal V.
\end{equation*}
\end{enumerate}
\end{lemma}


Clearly, (\ref{alg}) can be rewritten as
\begin{equation}\label{cons_comp}
x_i(k+1)=\sum_{j\in\mathcal V}a_{ij}(k) x_j(k)+\omega_i(k),
\end{equation}
where $\omega_i(k)=P_X(v_i(k)-\alpha_k d_i(k))-v_i(k)$.
In fact, (\ref{cons_comp}) is a consensus dynamics with disturbance  $\omega_i(k)$. Combining (\ref{distur}) with $\lim_{k\to\infty}\alpha_k=0$, we obtain $\lim_{t\to\infty} \omega_i(k)=0$. For such a dynamics, consensus can be achieved under a UJSC graph sequence as discussed in \cite{shi2013robust,wang2008robust}.
However, Lemma \ref{cons:lem} (ii), indicating the consensus rate, has not been proved in \cite{shi2013robust,wang2008robust}.

Referring to Theorem 1 in \cite{lin2016distributed}, we have the following result for (\ref{alg}) under a fixed digraph.

\begin{lemma}\label{con:fix}
Consider the graph sequence given by
$\mathcal G(k)=\mathcal G_A$
for $k\in \mathbb N$,
where $\mathcal G_A$ is a strongly connected  graph associated with weight matrix $A=[a_{ij}]$.
Suppose that $\sum_{j\in\mathcal V}a_{ij}=1$ for all $i \in \mathcal V$.
Under Assumptions \ref{con_ass}, \ref{lip_ass} and \ref{step_ass}, algorithm (\ref{alg}) achieves a solution to 
\begin{equation}\label{pro:fixed}
\min~~\sum_{i\in\mathcal V} \mu_i(A) f_i(x),
\quad {\rm s.t.}~~x\in X,
\end{equation}
where $\mu(A)=[\mu_1(A),\dots, \mu_n(A)]'$ is the Perron vector of $A$ such that $\mu(A)'A=\mu(A)'$.
\end{lemma}

\begin{remark}
Lemma \ref{con:fix} implies that (\ref{alg}) optimizes a weighted sum of the local cost functions if $\mathcal G(k)$ is fixed.
In fact, the result can be directly extended as follows.
Let $\mathcal G_{A_1}, \dots, \mathcal G_{A_p}$
be strongly connected graphs with an identical Perron vector $\mu(A)$.
Consider a time-varying graph sequence $\mathcal G(k)$, which switches within $\{\mathcal G_{A_1},\dots,\mathcal G_{A_p}\}$, i.e., $\mathcal G(k)\in \{\mathcal G_{A_1},\dots,\mathcal G_{A_p}\}$ for all $k \in \mathbb N$.
If Assumptions \ref{con_ass}, \ref{lip_ass}, \ref{wei_ass} and \ref{step_ass} hold, 
then algorithm (\ref{alg}) achieves a solution to (\ref{pro:fixed}). 
The proof is similar to that of Theorem 1 in \cite{lin2016distributed},
and is omitted here.
\end{remark}

\subsection{Convergence analysis}
In this section, we analyze whether (\ref{alg}) is still convergent in the absence of doubly stochastic weight matrices for the communication network. 

Let $\mathcal G_{A_1},\dots, \mathcal G_{A_p}$ be strongly connected graphs with weight matrices $A_1,\dots, A_p$, respectively. 
We define $\mu(A_l)=[\mu_1(A_l),\dots, \mu_n(A_l)]'$ as the Perron vectors of $A_l$ for $l\in \{1,\dots,p\}$, and moreover,
\begin{equation}\label{X_mu}
X_{\mu(A_l)}^* \triangleq 
\Big\{z~|~z={\rm argmin}_{x\in X} 
\sum_{i\in\mathcal V}\mu_i(A_l) f_i(x)\Big\}.
\end{equation}
Then we have the following result, whose proof can be found in Section \ref{pf:non_con}.

\begin{theorem}\label{non_con_the}
Let Assumptions \ref{con_ass}, \ref{lip_ass}, \ref{wei_ass} and \ref{step_ass} hold.
Consider a time-varying graph $\mathcal G(k)$, which switches within $\{\mathcal G_{A_1},\dots, \mathcal G_{A_p}\}$, i.e., $\mathcal G(k) \in \{\mathcal G_{A_1},\dots, \mathcal G_{A_p}\}$ for all $k \in \mathbb N$.
If $\bigcap_{l=1}^p X_{\mu(A_l)}^*=\emptyset$,
then there exists a graph sequence 
$\{\mathcal G(k)\}$ such that
algorithm (\ref{alg}) is not convergent.
\end{theorem}

It follows from Lemma \ref{con:fix} that (\ref{alg}) converges to $X_{\mu(A_l)}^*$ if $\mathcal G(k) =\mathcal G_{A_l}$ for all $k \in \mathbb N$.
As a result, $x_i(k)$ intends to oscillate
if $\mathcal G(k)$ switches within $\{\mathcal G_{A_1}, \dots, \mathcal G_{A_p}\}$. 
This leads to the non-convergent result.

\subsection{Condition for convergence}
Inspired by the above observation, we explore a condition to guarantee the convergence of (\ref{alg}) here. We present the main result in the following theorem, whose proof is given in Section \ref{pf:nec_suf}.

\begin{theorem}\label{con:nec_suf}
Let Assumptions \ref{con_ass}-\ref{step_ass} hold. Then algorithm (\ref{alg}) always achieves a solution to (\ref{pro_for}) for any UJSC graph sequence $\{\mathcal G(k)\}$ if and only if 
$\bigcap_{i\in \mathcal V} X_i^*\not= \emptyset$, where 
\begin{equation}
X_i^*=\big\{z~|~z={\rm argmin}_{x\in X} f_i(x)\big\}.
\end{equation}
\end{theorem}

\begin{remark}
Consider $\mathcal G(k)$ switching within
$\{\mathcal G_{A_1},\dots, \mathcal G_{A_p}\}$, where $\mathcal G_{A_l}$ is a
strongly connected graph for $l\in\{1,\dots,p\}$.
From the proof of Theorem \ref{con:nec_suf}, we conclude that if  weight matrix $A_l$ associated with $\mathcal G_{A_l}$ can be chosen freely under Assumptions \ref{gra_ass} and \ref{wei_ass}, then $\bigcap_{l=1}^p X_{\mu(A_l)}^*=\emptyset$ 
if and only if
$\bigcap_{i \in \mathcal V}X_i^*=\emptyset$.
This implies that the Assumption $\bigcap_{i \in \mathcal V}X_i^*=\emptyset$ in Theorem \ref{non_con_the} can be cast into
$\bigcap_{i \in \mathcal V}X_i^*=\emptyset$ 
in some sense.
\end{remark}

\begin{remark}
Note that Theorem \ref{con:nec_suf} provides a necessary and sufficient condition for the convergence of (\ref{alg}). 
The necessity can be inferred from Theorem \ref{non_con_the}, while the sufficiency is also considerable.
To be specific, if $\bigcap_{i\in \mathcal V} X_i^*\not= \emptyset$, all agents reach a consensus solution in $\bigcap_{i\in \mathcal V} X_i^*$
because agent $i$ intends to achieve consensus with its neighbors, and meanwhile, forces its state $x_i(k)$ to be in $X_i^*$.
\end{remark}

It is worthwhile to mention that Theorem \ref{con:nec_suf} is an extension of results for the convex intersection computation problem \cite{lou2014approximate,shi2012reaching} as follows.
In a network of $n$ nodes, all agents attempt to seek a consensus point in $\bigcap_{i\in\mathcal V} \Omega_i$ distributedly, where agent $i$ only knows its local convex set $\Omega_i$, and $\bigcap_{i\in\mathcal V}\Omega_i \not=\emptyset$.
Suppose that $\sum_{k=0}^\infty \alpha_k=\infty$. By the algorithms proposed in \cite{lou2014approximate, shi2012reaching}, the goal is achieved under Assumptions \ref{gra_ass} and \ref{wei_ass}.
Define 
$f_i(x)={\rm dist}^2(x, \Omega_i)$ and 
$X\in \mathbb R^m$.
Then the convex intersection computation problem is a special case of (\ref{pro_for}).
In fact, $\bigcap_{i\in \mathcal V} X_i^*\not= \emptyset$ means $\bigcap_{i\in\mathcal V} \Omega_i\not= \emptyset$.
Therefore, the sufficiency in Theorem \ref{con:nec_suf} is an extension for convergence results shown in \cite{lou2014approximate,shi2012reaching}. In addition, the necessity has not been proved in \cite{lou2014approximate,shi2012reaching}.

\subsection{Periodic graph sequences}
Theorems \ref{non_con_the} and \ref{con:nec_suf} indicate that the convergence of (\ref{alg}) cannot be guaranteed in general if the graph sequence $\{\mathcal G(k)\}$ can be chosen and switched freely. 
However, whether is it still convergent for some special graph sequences?
In this subsection, we investigate properties of (\ref{alg}) under periodic and quasi-periodic graph sequences.

Let $\mathcal G_{A_l}$ be a graph associated with weight matrix $A_l$ for $l\in \{1,\dots,p\}$, where $p\ge 2$.
Consider $\mathcal G(k)$ switching periodically
within $\{\mathcal G_{A_1},\dots, \mathcal G_{A_p}\}$.
To be specific, the graph sequence is given by
$\mathcal G(tp+l-1)=\mathcal G(l-1)=\mathcal G_{A_l}$ for all $t\in \mathbb N$.
For simplification, we write the sequence as
$$\mathcal G_{A_1} \rightarrow \cdots \rightarrow \mathcal G_{A_p} \rightarrow \mathcal G_{A_1} \rightarrow \cdots \rightarrow \mathcal G_{A_p} \rightarrow \cdots.$$
Let $\mu^1,\mu^2,\dots, \mu^p$ be Perron vectors of 
$$(A_pA_{p-1}\cdots A_1),~ 
(A_{p-1}A_{p-2}\cdots A_1A_p), ~\dots, ~(A_1A_pA_{p-1}\cdots A_2),$$ respectively.
Then we have the following result for  (\ref{alg}), whose proof is provided in Section \ref{pf:con_per}.

\begin{theorem}\label{con:periodic}
Let Assumptions \ref{con_ass}-\ref{step_ass} hold.
If the set constraint $X$ is compact and the stepsize sequence $\{\alpha_k\}$ is non-increasing,
then algorithm (\ref{alg}) achieves a solution to 
\begin{equation}\label{per_for}
\min~~\sum_{i\in\mathcal V}
\frac 1p(\mu^1_i+\dots+\mu^p_i) f_i(x),~~~~
{\rm s.t.}~~x\in X,
\end{equation}
where $\mu^l_i$ is the $i$-th entry of $\mu^l$.
\end{theorem}

\begin{remark}
Theorem \ref{con:periodic} indicates that (\ref{alg}) is convergent under a periodically switching general graph sequence $\{\mathcal G(k)\}$, and moreover, the converged solution minimizes a weighted sum of local cost functions, where the weights depend on the Perron vectors of some product matrices of the underlying periodically switching  graphs. 
\end{remark}

Intuitively, (\ref{alg}) is convergent under a periodic graph sequence because at each time interval $[tp,(t+1)p)$, the joint graph 
$\mathcal G([tp,(t+1)p))$ is time-invariant. Recalling Lemma \ref{con:fix} gives the convergence result. 
However, it should be note that Theorem \ref{con:periodic} is not so straightforward.
The intuition indicates that (\ref{alg}) may achieve a solution to  
\begin{equation}\label{intu_res}
\min~~\sum_{i\in\mathcal V}
{\hat\mu}_i f_i(x),
\quad {\rm s.t.}~~x\in X,
\end{equation}
where $\hat\mu=[\hat\mu_1,\dots,\hat\mu_n]'$ is the Perron vector of $(A_1A_2\cdots A_p)$.
Clearly, the solution set to (\ref{intu_res}) is generally different from that of (\ref{per_for}), which makes a contradiction. Take $p=2$ for interpretation.
Consider $f_i(x)$ being strictly convex. Then there is a unique solution to any weighted sum of the local cost functions. 
If (\ref{alg}) achieves a solution to a weighted optimization problem, $\lim_{k\to\infty}x_i(k)$ is independent of the initial state. 
As a result, (\ref{alg}) reaches the same solution under both graph sequences
$$\mathcal G_{A_1} \rightarrow \mathcal G_{A_2} \rightarrow \mathcal G_{A_1}\rightarrow \mathcal G_{A_2} \rightarrow\cdots~~ {\rm and} ~~
\mathcal G_{A_2} \rightarrow \mathcal G_{A_1}\rightarrow \mathcal G_{A_2} \rightarrow \mathcal G_{A_1} \rightarrow \cdots.$$
However, the Perron vectors of $A_1A_2$ and $A_2A_1$ are generally not identical,
which leads to (\ref{intu_res}) with different solutions under the above two sequences.
This implies the incorrectness of (\ref{intu_res}).
By the proof of Theorem \ref{con:periodic}, we conclude that $\sum_{l=1}^p \mu^l$ is independent of the initial point of the graph sequence, which verifies the correctness of Theorem \ref{con:periodic}.

In \cite{lobel2010distributed}, the authors considered the weight matrix $A(k)$ drawn independently from a probability space, and explored convergence properties of a distributed subgradient algorithm.
It was shown that the convergence relied on the expected graph.
Inspired by the result, it seem that under periodic graphs,
(\ref{alg}) may achieve a solution to
\begin{equation}\label{exp:con_pro}
\min~~\sum_{i\in\mathcal V}
\mu^e_i f_i(x),
\quad {\rm s.t.}~~x\in X,
\end{equation}
where $\mu^{e}=[\mu^e_1,\dots,\mu^e_n]'$ is the Perron vector of $\sum_{l=1}^p \frac 1p A_l$.
However, Theorem \ref{con:periodic} demonstrates that the convergence of (\ref{alg}) under periodic graph sequences is different from that shown in \cite{lobel2010distributed}.

It follows from Theorem \ref{con:periodic} that 
a periodic graph sequence is a sufficient condition to guarantee the convergence of (\ref{alg}).
It is natural to consider whether the condition is also necessary.
We relax the periodic condition slightly, and define a broader class of quasi-periodic graph sequences as follows.
Let $\mathcal G_{A_l}$ be a graph associated with weight matrix $A_l$ for $l\in \{1,\dots,p\}$, where $p\ge 2$.
$\{\mathcal G(k)\}$ is called a quasi-periodic sequence if
it switches within $\{\mathcal G_{A_1},\dots, \mathcal G_{A_p}\}$ at each time interval 
$[tp, (t+1)p)$ for $t\in \mathbb N$,
but the order of $\mathcal G_{A_l}$ can be changed over $t$.
For instance, we consider $p=3$. Then we can take the graph sequence as  $\mathcal G_{A_1} \rightarrow \mathcal G_{A_2}\rightarrow \mathcal G_{A_3}$ at the time interval $[0, 3)$, and $\mathcal G_{A_1} \rightarrow \mathcal G_{A_3}\rightarrow \mathcal G_{A_2}$ at the time interval $[3, 6)$.
The following theorem, proved in Section \ref{pf:noper}, addresses a property of (\ref{alg}) under quasi-periodic graph sequences.

\begin{theorem}\label{con:nonperiodic}
Let Assumptions \ref{con_ass}-\ref{step_ass} hold.
Suppose that the set constraint $X$ is compact, the stepsize sequence $\{\alpha_k\}$ is non-increasing, and moreover, $\bigcap_{i\in \mathcal V} X_i^*=\emptyset$.
Then algorithm (\ref{alg}) is convergent for any quasi-periodic graph sequence if and only if $p=2$.
\end{theorem}

\begin{remark}
By Theorem \ref{con:periodic}, 
the optimization problem (\ref{per_for})
relies on the order of the graph sequence at each time interval $[tp, (t+1)p)$ if $p\ge 3$. With the help of Theorem \ref{non_con_the}, Theorem \ref{con:nonperiodic} can be inferred.
\end{remark}

Here, we relax the periodic graph by another way. 
Let $\mathcal G_{A_l}$ be a graph associated with weight matrix $A_l$ for $l\in \{1,\dots,p\}$, where $p\ge 2$.
Consider $\mathcal G(k)$ switching within 
$\{\mathcal G_{A_1},\dots, \mathcal G_{A_p}\}$ 
at each time interval $[tD, (t+1)D)$, where $D> p$.
However, $\mathcal G_{A_l}$ may appear with different frequencies at time intervals $[tD, (t+1)D)$ for $t\in \mathbb N$.
For instance, we consider $D=3$ and $p=2$. We can take the graph sequence as $\mathcal G_{A_1} \rightarrow \mathcal G_{A_1}\rightarrow \mathcal G_{A_2}$ at the time interval $[0, 3)$, and $\mathcal G_{A_1} \rightarrow \mathcal G_{A_2}\rightarrow \mathcal G_{A_2}$ at the time interval $[3, 6)$. 
In this case, we have the following corollary.

\begin{corollary}\label{per_coro}
Let Assumptions \ref{con_ass}-\ref{step_ass} hold. Suppose that $\bigcap_{i\in \mathcal V} X_i^*=\emptyset$.
If $A_l$  can be chosen freely, then there exists a graph sequence such that algorithm (\ref{alg}) is not convergent.
\end{corollary}

\emph{Proof.} The main idea on the proof  focuses on discussing the Perron vector at each time interval $[tD, (t+1)D)$. It is similar to that of the case of $p\ge 3$ in Theorem \ref{con:nonperiodic}, and is omitted here.

\begin{remark}
Theorem \ref{con:nonperiodic} and Corollary \ref{per_coro} indicate that (\ref{alg}) is not convergent in general if a periodic graph sequence is with a slight modification.
Therefore, we can conclude that the periodic condition is very important to guarantee the convergence of (\ref{alg}) for general graph sequences.
\end{remark}

\section{Proofs}
In this section, we introduce several useful lemmas, and then prove the results presented in the last section.

\subsection{Supporting lemmas}
Referring to \cite{lou2015nash,nedic2009distributed}, we define transition matrices as
\begin{equation*}
\Phi(k,s)=A(k)A(k-1)\cdots A(s)
\end{equation*}
for $s,k\in \mathbb N$ with $k\ge s$, where $\Phi(k,k)=A(k)$. Then $\Phi(k,s)$ is a stochastic matrix. In light of Lemma 2 in \cite{nedic2009distributed}, we have the following result.

\begin{lemma}\label{tran_mat_pro}
Under Assumptions \ref{gra_ass} and \ref{wei_ass}, $[\Phi(s+(n-1)B-1,s)]_{ij} \ge \eta^{(n-1)B}$ for all $i,j\in \mathcal V$ and $s\in \mathbb N$.
\end{lemma}


The following lemma, found from Lemma 3 in \cite{hajnal1958weak}, will be used for the consensus analysis.
\begin{lemma}\label{cons_lem}
For $\mu=[\mu_1,\dots, \mu_n]'\in \mathbb R^n$, we define $g(\mu)=\max_{1\le i,j\le n}|\mu_i-\mu_j|$.
If $P=[p_{ij}]\in \mathbb R^{n\times n}$ 
is a stochastic matrix, then 
$g(P\mu) \le \tau(P)g(\mu)$,
where $\tau(P)=1-\min_{i,j}\sum_{s=1}^n 
\min\{p_{is},p_{js}\}$.
\end{lemma}

We introduce two lemmas about infinite series for the convergence analysis. The first one is a deterministic version of Lemma $11$ on page $50$ in \cite{polyak1987introduction}, while the second one is collected from Lemma 7 in \cite{nedic2010constrained}.
\begin{lemma}\label{sequ_lem}
Let $\{a_k\}, \{b_k\}$ and $\{c_k\}$
be non-negative sequences with 
$\sum_{k=0}^\infty b_k<\infty$. 
If $a_{k+1}\le a_k+b_k-c_k$
holds for all $k\in \mathbb N$, 
then the limit $\lim_{k\to \infty}a_k$ exists and
is a finite number.
\end{lemma}

\begin{lemma}\label{ser_lem}
Let $0<\beta< 1$ and $\{\gamma_k\}$ be a positive scalar sequence. 
\begin{enumerate}[(i)]
\item If $\lim_{k\to \infty} \gamma_k=0$, then
      $\lim_{k\to\infty}\sum_{s=0}^k \beta^{k-s}\gamma_s=0$. 
\item If $\sum_{k=0}^\infty \gamma_k  
      <\infty$, then
      $\sum_{k=0}^{\infty}\sum_{s=0}^k \beta^{k-s}\gamma_s <\infty$.
\end{enumerate}
\end{lemma}

Similar to Lemma 6 in \cite{nedic2009distributed}, we have the following result for (\ref{alg}) under general graphs. Since it will be frequently used later, we provide a concise proof here. 
\begin{lemma}\label{ite_lem}
Let $x_i(k)$ be generated by algorithm (\ref{alg}). 
Suppose that Assumptions \ref{con_ass}, \ref{lip_ass} and \ref{wei_ass} (i) hold.
For $z\in X$ and $k\in \mathbb N$, we have
\begin{equation}\label{ite_ine}
|x_i(k+1)-z|^2 \le 
\sum_{j\in\mathcal V}a_{ij}(k)|x_j(k)-z|^2
+\alpha_k^2 L^2
-2\alpha_k(f_i(v_i(k))-f_i(z)).
\end{equation}
\end{lemma}

\emph{Proof.}
Clearly, (\ref{alg}) can be rewritten as
\begin{equation}\label{com_alg3}
x_i(k+1)=v_i(k)-\alpha_k d_i(k)+\varphi_i(k),
\end{equation}
where $\varphi_i(k)=P_X (v_i(k)-\alpha_k d_i(k))-(v_i(k)-\alpha_k d_i(k))$.
Recalling (\ref{pro_ine2}) gives
\begin{equation}\label{lem_ine}
\begin{aligned}
|x_i(k+1)-z|^2 
\le&|v_i(k)-\alpha_k d_i(k)-z|^2
-|\varphi_i(k)|^2 \\ 
\le&|v_i(k)-z|^2+\alpha_k^2|d_i(k)|^2 -2\alpha_k(f_i(v_i(k))-f_i(z))
-|\varphi_i(k)|^2.
\end{aligned}
\end{equation}

Because the norm square function is convex, 
$\sum_{j\in \mathcal V} a_{ij}(k) |x_j(k)-z|^2 \ge |v_i(k)-z|^2$. 
By combining (\ref{lip_ine}) and (\ref{lem_ine}), the conclusion follows.
$\hfill\Box$

\subsection{Proof of Lemma \ref{cons:lem}}
\label{fix:con_lem_pf}
Here, we consider $m=1$ to simplify the proof, and otherwise, we can use Kronecker product when necessary. 
By (\ref{cons_comp}), we have
\begin{equation}\label{com_alg}
x(k+1)=A(k)x(k)+\omega(k)
=\Phi(k,s)x(s)+
\sum_{r=s}^{k-1}\Phi(k,r+1)\omega(r)
+\omega(k),
\end{equation}
where $x(k)=[x_1(k),\dots, x_n(k)]'$, and  $\omega(k)=[\omega_1(k),\dots, \omega_n(k)]'$.
Recalling (\ref{pro_ine1}) yields
\begin{equation}\label{distur}
\begin{aligned}
|\omega_i(k)| \le 
|v_i(k) -\alpha_k d_i(k)-v_i(k)| \le
\alpha_k L.
\end{aligned}
\end{equation}

Take $h(k)=\max_{i,j\in\mathcal V}|x_i(k)-x_j(k)|$ and $T=(n-1)B$. 
It follows from (\ref{com_alg}), Lemma \ref{cons_lem} and 
$g(\mu+\nu)\le g(\mu)+2\max_i |\nu_i|$ that
\begin{equation}\label{con:ine1}
h(s+(t+1)T) \le \tau\big(\Phi(s+(t+1)T-1, s+tT)\big)h(s+tT)+
\sum_{r=s+tT}^{s+(t+1)T-1}2L\alpha_r.
\end{equation}
In light of Lemma \ref{tran_mat_pro},
$\Phi(s+(t+1)T-1, s+tT)\ge \eta^T$, 
and then $\tau\big(\Phi(s+(t+1)T-1, s+tT)\big) \le 1-\eta^T$.
Define $\beta_{s,t}=\sum_{r=s+tT}^{s+(t+1)T-1}\alpha_r$. By (\ref{con:ine1}), we obtain
\begin{equation}\label{con:ine2}
\begin{aligned}
h(s+(t+1)T) \le& (1-\eta^T)h(s+tT)+
2L\beta_{s,t} \\
\le& (1-\eta^T)^{t+1} h(s)
+\sum_{r=0}^{t} 2L(1-\eta^T)^{t-r} \beta_{s,r}.
\end{aligned}
\end{equation}

If $\lim_{t\to\infty}\alpha_t=0$, then
$\lim_{t\to\infty}\beta_{s,t}=0$.
Due to Lemma \ref{ser_lem} (i),  $\lim_{t\to\infty}\sum_{r=0}^t 
(1-\eta^T)^{t-r} \beta_{s,r}=0$.
Clearly, $\lim_{t\to\infty}(1-\eta^T)^{t+1} h(s)=0$.
Thus, $\lim_{k\to\infty}h(k)
=0$. By the definition of $y(k)$,  $h(k)\ge|x_i(k)-y(k)|$, and then $\lim_{k\to\infty}|x_i(k)-y(k)|=0$.

Combining (\ref{con:ine2}) with $2\alpha_{s+tT}\beta_{s,r} \le \alpha_{s+tT}^2 +\beta_{s,r}^2$,
we derive
\begin{equation}
\begin{aligned}
\sum_{k=0}^{\infty}\alpha_kh(k)&=
\sum_{t=0}^{\infty}\sum_{s=0}^{T-1}
\alpha_{s+tT}h(s+tT)\\
&\le \sum_{s=0}^{T-1}\sum_{t=0}^\infty
\Big[\alpha_{s+tT}
(1-\eta^T)^t h(s)
+\sum_{r=0}^{t-1}
(1-\eta^T)^{(t-1)-r}\beta_{s,r}^2 L\\
&\qquad+\alpha_{s+tT}^2\sum_{r=0}^{t-1}
(1-\eta^T)^{(t-1)-r}L
\Big].
\end{aligned}
\end{equation}

Because of $\sum_{k=0}^\infty \alpha_k^2 <\infty$, $\alpha_k$ is bounded,
and then $\sum_{t=0}^\infty \alpha_{s+tT}(1-\eta^T)^t h(s)$ 
is also bounded.
Obviously, $\beta_{s,t}^2 \le \sum_{r=s+tT}^{s+(t+1)T-1}2\alpha_r^2$.
By $\sum_{k=0}^\infty \alpha_k^2 <\infty$, we have $\sum_{t=0}^{\infty}\beta_{s,t}^2 <\infty$.
Recalling \ref{ser_lem} (ii) gives
$\sum_{t=0}^\infty\sum_{r=0}^{t-1} (1-\eta^T)^{(t-1)-r}\beta_{s,r}^2L<\infty$.
In addition, $\sum_{r=0}^{t-1}
(1-\eta^T)^{(t-1)-r}L\le L/\eta^T$, and then
$\sum_{t=0}^{\infty}\alpha_{s+tT}^2 \sum_{r=0}^{t-1}(1-\eta^T)^{(t-1)-r}L <\infty$.
Therefore, 
$\sum_{k=0}^{\infty} \alpha_k h(k)< \infty$, and moreover,
$\sum_{k=0}^\infty \alpha_k|x_i(k)-y(k)|<\infty$.
This completes the proof.
$\hfill\Box$

\subsection{Proof of Theorem \ref{non_con_the}}
\label{pf:non_con}
Because of the convexity of $f_i$ and $\mu_i(A_l)>0$, $\sum_{i\in\mathcal V}\mu_i(A_l) f_i(x)$ is a convex function, and as a result, $X_{\mu(A_l)}^*$ is a closed convex set.
Define $x(k)=[x_1'(k), \dots, x_n'(k)]'$.
It follows from Lemma \ref{con:fix} that (\ref{alg})
converges to a point in $X_{\mu(A_l)}^*$ under the fixed graph sequence 
$\mathcal G(k)=\mathcal G_{A_l}$ for $k\in\mathbb N$.
Thus, for any $\epsilon>0$ and initial point $x(0)\in \mathbb R^{mn}$, there must be 
$T_l(\epsilon, x(0))\in \mathbb N$ 
such that
\begin{equation*}
{\rm dist}(x_i(t), X_{\mu(A_l)}^*)<\epsilon, 
~\forall i\in \mathcal V,
~\forall t\ge T_l(\epsilon, x(0)).
\end{equation*}

Due to $\bigcap_{l=1}^p X_{\mu(A_l)}^* =\emptyset$, there must be
$i,j\in\{1,\dots,p\}$ and a scalar \\
$d \triangleq {\rm dist}(X_{\mu(A_i)}^*, X_{\mu(A_j)}^*)>0$.
Without loss of generality, we assume $i=1$ and $j=2$.
Then we construct  time sequences $\{t_k\}$ and $\{s_k\}$,
and a switching graph sequence 
$\{\mathcal G(k)\}$ as follows.

Let $s_0=0, t_0=s_0$ and $x(0)\in \mathbb R^{mn}$. Furthermore,
\begin{equation}\label{graph_seq_cons}
\begin{aligned}
s_1&=T_1(d/3, x(t_0)), t_1=t_0+s_1, 
{\rm and}~\mathcal G(k)=\mathcal G_{A_1} 
~{\rm for}~k=t_0+1,\dots, t_1;\\
s_2&=T_2(d/3, x(t_1)), t_2=t_1+s_2, 
{\rm and}~\mathcal G(k)=\mathcal G_{A_2}  
~{\rm for}~k=t_1+1,\dots, t_2;\\
&~~\vdots \\
s_{2k+1}&=T_1(d/3, x(t_{2k})), t_{2k+1}=t_{2k}+s_{2k+1}, \\
&~~~~{\rm and}~\mathcal G(k)=\mathcal G_{A_1} 
~{\rm for}~k=t_{2k}+1,\dots, t_{2k+1};\\
s_{2k+2}&=T_2(d/3, x(t_{2k+1})), t_{2k+2}=t_{2k+1}+s_{2k+2}, \\
&~~~~{\rm and}~\mathcal G(k)=\mathcal G_{A_2}
~{\rm for}~k=t_{2k+1}+1,\dots, t_{2k+2}.
\end{aligned}
\end{equation}
Then $|x(t_{2k+1})-x(t_{2k+2})|>d/3$ 
for all $k \in \mathbb N$.
Therefore, $x(t)$ is not convergent, and the conclusion holds.
$\hfill\Box$

\subsection{Proof of Theorem \ref{con:nec_suf}}
\label{pf:nec_suf}
\emph{(Necessity)}
The necessary is shown by contradiction. To be specific, if $\bigcap_{i\in \mathcal V}X_i^*=\emptyset$, there always exists a graph sequence such that (\ref{alg}) is not convergent.

By Lemma \ref{sto_mat_lem}, for any positive stochastic vector $\mu\in \mathbb R^n$, there exists a stochastic matrix $A$, whose Perron vector is $\mu$. Moreover, the graph $\mathcal G_A$ associated with $A$ is strongly connected.
Here, we take two positive stochastic vectors $\mu(A_1)$ and $\mu(A_2)$ associated with  matrices $A_1$ and $A_2$, and graphs $\mathcal G_{A_1}$ and $\mathcal G_{A_2}$ as follows.
Define $X_{\mu(A_1)}^*$ and $X_{\mu(A_2)}^*$ by (\ref{X_mu}).
Take $\mu(A_1)=[1/n,\dots,1/n]$.
In light of Lemma \ref{con:fix}, all agents converge to a point $\hat x$ in
$X_{\mu(A_1)}^*$  if 
$\mathcal G(k)=\mathcal G_{A_1}$. 
If $\bigcap_{i\in \mathcal V} X_i^*= \emptyset$,
then there must be $i_0\in \mathcal V$ and $\tilde x\in X_{i_0}^*$ such that
$f_{i_0}(\tilde x)< f_{i_0}({\hat x})$.
Take $\mu(A_2)$ such that
\begin{equation*}
\begin{aligned}
{\mu_{i_0}(A_2)}\big(f_{i_0}(\hat x)-f_{i_0}
(\tilde x)\big)
>\sum_{i\in\mathcal V, i\not= i_0} \mu_i(A_2) (f_i(\tilde x)-f_i(\hat x)),
\end{aligned}
\end{equation*}
where $\mu_i(A_2)$ is the $i$-th entry of $\mu(A_2)$. Consequently,
$\sum_{i\in\mathcal V}\mu_i(A_2)f_i(\tilde x)< \\
\sum_{i\in\mathcal V}\mu_i(A_2)f_i(\hat x)$.
Therefore, we conclude that $X_{\mu(A_1)}^*\bigcap X_{\mu(A_2)}^*=\emptyset$.
In view of Theorem \ref{non_con_the}, there exists a graph sequence such that (\ref{alg}) is not convergent.

\emph{(Sufficiency)} 
The sufficiency is proved by the following three steps.

\emph{Step 1.}
In this step, we show that $\{x_i(k)\}$ is bounded.
Define $X_s^*\triangleq\bigcap_{i\in\mathcal V}X_i^*\not= \emptyset$
and take $x^*\in X_s^*$. By setting $z=x^*$ in (\ref{ite_ine}), we derive
\begin{equation}\label{suff_ine1}
|x_i(k+1)-x^*|^2 \le 
\sum_{j\in\mathcal V}a_{ij}(k)|x_j(k)-x^*|^2
+\alpha_k^2 L^2
-2\alpha_k (f_i(v_i(k))-f_i(x^*)).
\end{equation}
Define $\xi(k)=\max_{i\in\mathcal V} |x_i(k)-x^*|^2$. Then
\begin{equation*}
\xi(k+1)\le \xi(k)+\alpha_k^2 L^2
-2\alpha_k \min\{f_i(v_i(k))-f_i(x^*)\}.
\end{equation*}
Notice that $f_i(v_i(k))-f_i(x^*) \ge 0$ and $\sum_{k=0}^\infty \alpha_k^2 L^2<\infty$. It follows from Lemma \ref{sequ_lem} that there exists $\xi^*$ such that 
$\lim_{k\to\infty} \xi(k)= \xi^*.$
As a result, $\{x_i(k)\}$ is bounded. 

Additionally, we have
\begin{equation*}
\max_{i\in\mathcal V}|x_i(k)-x^*|
-\max_{i,j\in\mathcal V}|x_i(k)-x_j(k)|
\le \min_{i\in\mathcal V}|x_i(k)-x^*|
\le \max_{i\in\mathcal V}|x_i(k)-x^*|.
\end{equation*}
Recalling $\lim_{k\to\infty}\max_{i,j\in\mathcal V}|x_i(k)-x_j(k)|=0$ gives
\begin{equation}\label{x_con_res}
\lim_{k\to\infty}|x_i(k)-x^*|=\xi^*,
~\forall i\in\mathcal V.
\end{equation}

Clearly, $y(k)$ is also bounded. 
In light of $\lim_{k\to \infty}|x_i(k)-y(k)|=0$, the sequence $\{|y(k)-x^*|\}$ is convergent for any $x^*\in X_s^*$.

\emph{Step 2.}
Define $\zeta_i(k)=|x_i(k)-x^*|^2$.
Recalling (\ref{suff_ine1}), we obtain
\begin{equation}\label{bou_ine1}
\begin{aligned}
\zeta_i&(k+1) \le
\sum_{j\in\mathcal V}[\Phi(k,s)]_{ij}\zeta_j(s)
+\sum_{r=s}^{k-1}\sum_{j\in\mathcal V}
[\Phi(k,r+1)]_{ij}\alpha_r^2 L^2 
+\alpha_k^2 L^2\\
&-\sum_{r=s}^{k-1}\sum_{j\in\mathcal V} 2[\Phi(k,r+1)]_{ij}\alpha_r
\big(f_j(v_j(r))-f_j(x^*)\big)-2\alpha_k\big(f_i(v_i(k))-f_i(x^*)\big).
\end{aligned}
\end{equation}

Clearly,
\begin{equation}\label{bou_ine2}
\begin{aligned}
-\sum_{r=s}^{k-1}&\sum_{j\in\mathcal V} [\Phi(k,r+1)]_{ij}\alpha_r
\big(f_j(v_j(r))-f_j(x^*)\big) \\
=&-\sum_{r=s}^{k-1}\sum_{j\in\mathcal V} [\Phi(k,r+1)]_{ij}\alpha_r
\big(f_j(v_j(r))-f_j(y(r))\big)\\
&-\sum_{r=s}^{k-1}\sum_{j\in\mathcal V} [\Phi(k,r+1)]_{ij}\alpha_r
\big(f_j(y(r))-f_j(x^*)\big).
\end{aligned}
\end{equation}

By (\ref{lip_ine}), $|f_j(v_j(r))-f_j(y(r))|<L|v_j(r)-y(r)|$.
In view of Lemma \ref{cons:lem} (ii), we obtain
$$\Big|\sum_{r=s}^{k-1}\sum_{j\in\mathcal V} [\Phi(k,r+1)]_{ij}\alpha_r
\big(f_j(v_j(r))-f_j(y(r))\big)\Big|<\infty.$$

\emph{Step 3.}
Here, we show that
$x_i(k)$ converges to a point in 
$X_s^*$ 
by contradiction. 
For any $x\in X$ and $\epsilon>0$ such that ${\rm dist}(x,X_i^*)>\epsilon$,
there must be $\delta>0$ such that $f_i(x)-f_i(x^*)>\delta$ due to the convexity and continuity of $f_i$.
By Lemma \ref{tran_mat_pro},
$\Phi(k,s)\ge \eta^{(n-1)B}$
for all $k\ge s+(n-1)B-1$.
For any $\epsilon>0$, we suppose
${\rm dist}(y(r),X_j^*)>\epsilon$.
Then for $k\ge s+(n-1)B$, we have
\begin{equation}\label{bou_ine4}
\begin{aligned}
-\sum_{r=s}^{k-1}&\sum_{j\in\mathcal V} [\Phi(k,r+1)]_{ij}\alpha_r
\big(f_j(y(r))-f_j(x^*)\big) \\
\le &-\sum_{r=s}^{k-(n-1)B}\sum_{j\in\mathcal V} [\Phi(k,r+1)]_{ij}\alpha_r
\big(f_j(y(r))-f_j(x^*)\big)\\
\le& -\delta\eta^{(n-1)B} \sum_{r=s}^{k-(n-1)B}\alpha_r.
\end{aligned}
\end{equation}

Substituting (\ref{bou_ine2}) and (\ref{bou_ine4}) into (\ref{bou_ine1}),
we obtain
$\lim_{k\to\infty}\zeta(k)=-\infty$ by Assumption \ref{step_ass}.
This contradicts with the boundedness of $x_i(k)$ proved in Step 1. 
Thus, \\
$\lim_{k\to\infty}\inf {\rm dist}(y(k), X_s^*)=0$.

Since $y(k)$ is bounded, it must have at least one limit point. In view of \\
$\lim_{k\to\infty}\inf {\rm dist}(y(k),X_s^*)=0$, one of the limit points, denoted by $y^*$, must be in $X_s^*$. As shown in Step 1, the sequence 
$\{|y(k)-y^*|\}$ is convergent, and as a result, the limit point is unique, i.e. $\lim_{k\to\infty}y(k)=y^*$.
By $\lim_{k\to\infty}|x_i(k)-y(k)|=0$, we conclude that for all $i\in\mathcal V$, $x_i(k)$ converges to the same $y^*\in X_s^*$.
This completes the proof.
$\hfill\Box$

\subsection{Proof of Theorem \ref{con:periodic}}
\label{pf:con_per}
Consider $p=2$ for simplification, and note that the idea can be directly extended to the case of $p>2$. We divide the proof into three steps as follows.

\emph{Step 1.}
Without loss of generality, we consider $\mathcal G(2k)=\mathcal G_{A_1}$ and
$\mathcal G(2k+1)=\mathcal G_{A_2}$ for  $k\in \mathbb N$.
Recalling (\ref{ite_ine}) gives
\begin{equation}\label{k1_ite_ine}
|x_i(2k+1)-z|^2 \le 
\sum_{j\in\mathcal V}[A_1]_{ij}|x_j(2k)-z|^2
+\alpha_{2k}^2 L^2
-2\alpha_{2k}\big(f_i(v_i(2k))-f_i(z)\big),
\end{equation}
and moreover,
\begin{equation}\label{k2_ite_ine}
\begin{aligned}
|x_i(2k+2)-z|^2 \le 
\sum_{j\in\mathcal V} [A_2]_{ij}|x_j(2k&+1)-z|^2
+\alpha_{2k+1}^2 L^2\\
&-2\alpha_{2k+1}
\big(f_i(v_i(2k+1))-f_i(z)\big).
\end{aligned}
\end{equation}
Notice that $\sum_{j\in \mathcal V} [A_2]_{ij}=1$. Substituting (\ref{k1_ite_ine}) into (\ref{k2_ite_ine}), we obtain
\begin{equation}\label{per_ite_ine}
\begin{aligned}
|x_i(2k&+2)-z|^2 \le
\sum_{j\in\mathcal V} [A_2A_1]_{ij}|x_j(2k)-z|^2
+(\alpha_{2k}^2+\alpha_{2k+1}^2 ) L^2\\
&-\sum_{j\in\mathcal V} 2\alpha_{2k}[A_2]_{ij}(f_j(v_j(2k))-f_j(z))
-2\alpha_{2k+1}(f_i(v_i(2k+1))-f_i(z)).
\end{aligned}
\end{equation}

Let $\mu^1=[\mu^1_1,\dots, \mu^1_n]'$, 
$\mu^2=[\mu^2_1,\dots, \mu^2_n]'$ 
be the Perron vectors of 
$A_2A_1$ and $A_1A_2$ such that 
$(\mu^1)'A_2A_1=(\mu^1)'$ and 
$(\mu^2)'A_1A_2=(\mu^2)'$, respectively. 
As a result, we have
\begin{equation*}
[(\mu^1)'A_2](A_1A_2)=[(\mu^1)'A_2],
~{\rm and}~
[(\mu^2)'A_1](A_2A_1)=[(\mu^2)'A_1].
\end{equation*}
Therefore, $(\mu^1)'A_2$, $(\mu^2)'A_1$ are the Perron vectors of $A_1A_2$ and $A_2A_1$, respectively.
Because the joint graph 
$\mathcal G_{A_1}\cup \mathcal G_{A_2}$
is strongly connected, 
the Perron vectors of both $A_2A_1$ and $A_1A_2$ are unique by the Perron-Frobenius theorem.
Thus, 
\begin{equation}
(\mu^2)'=(\mu^1)'A_2, ~~{\rm and}~~ (\mu^1)'=(\mu^2)'A_1.
\end{equation}

Define $X_p^*=\{z~|~z={\rm argmin}_{x\in X} \sum_{i\in\mathcal V}\frac 12 (\mu^1_i+\mu^2_i)f_i(x)\}$. 
Let $x^*\in X_p^*$ and take $z=x^*$.
Multiplying $\mu^1_i$ to both sides of 
(\ref{per_ite_ine}) and summing all $i\in \mathcal V$, we obtain
\begin{equation*}\label{per_ite1}
\begin{aligned}
\sum_{i\in\mathcal V}& \mu^1_i|x_i(2k+2)-x^*|^2 \le
\sum_{i\in\mathcal V} \mu^1_i|x_i(2k)-x^*|^2
+\sum_{i\in\mathcal V} \mu^1_i
(\alpha_{2k}^2+\alpha_{2k+1}^2) L^2\\
&-\sum_{i\in\mathcal V}2  
\mu^2_i\alpha_{2k}
\big(f_i(v_i(2k))-f_i(x^*)\big)
-\sum_{i\in\mathcal V} 
2 \mu^1_i \alpha_{2k+1}
\big(f_i(v_i(2k+1))-f_i(x^*)\big).
\end{aligned}
\end{equation*}
By a similar procedure for discussing $|x_i(2k+2)-x^*|$, we also have
\begin{equation*}\label{per_ite2}
\begin{aligned}
\sum_{i\in\mathcal V} \mu^2_i|x_i(2k+3)-x^*|^2 \le
\sum_{i\in\mathcal V} \mu^2_i|x_i(2k+1)-x^*|^2
+\sum_{i\in\mathcal V} \mu^2_i
(\alpha_{2k+1}^2+\alpha_{2k+2}^2)L^2&\\
-\sum_{i\in\mathcal V} 
2\mu^1_i\alpha_{2k+1}
\big(f_i(v_i(2k+1))-f_i(x^*)\big)
-\sum_{i\in\mathcal V}
2\mu^2_i\alpha_{2k+2}
\big(f_i(v_i(2k+2))-f_i(x^*&)\big).
\end{aligned}
\end{equation*}

Define $\chi_k=\sum_{i\in\mathcal V}
\mu^1_i|x_i(2k)-x^*|^2 
+\sum_{i\in\mathcal V} \mu^2_i|x_i(2k+1)-x^*|^2$. 
Notice that $\sum_{i\in \mathcal V}\mu^1_i=1$ and $\sum_{i\in \mathcal V}\mu^2_i=1$. Combining the above two inequalities, we derive 
\begin{equation}\label{per_ine1}
\begin{aligned}
\chi_{k+1} \le &\chi_k
+(\alpha_{2k}^2+2\alpha_{2k+1}^2
+\alpha_{2k+2}^2) L^2
-\sum_{i\in\mathcal V}2 (\mu^1_i+\mu^2_i)\alpha_{2k}
\big(f_i(v_i(2k))-f_i(x^*)\big) \\
&-\sum_{i\in\mathcal V}
2(\mu^1_i+\mu^2_i)\alpha_{2k+1}
\big(f_i(v_i(2k+1))-f_i(x^*)\big)
+M_1(k)+M_2(k),
\end{aligned}
\end{equation}
where 
$$M_1(k)=\sum_{i\in\mathcal V}2
\mu^1_i\big[\alpha_{2k}
\big(f_i(v_i(2k))-f_i(x^*)\big)
-\alpha_{2k+1}
\big(f_i(v_i(2k+1))-f_i(x^*)\big)\big],$$
and moreover,
$$M_2(k)=\sum_{i\in\mathcal V} 
2\mu^2_i\big[ \alpha_{2k+1}
\big(f_i(v_i(2k+1))-f_i(x^*)\big)
-\alpha_{2k+2}
\big(f_i(v_i(2k+2))-f_i(x^*)\big)\big].$$

\emph{Step 2.} In this step, we analyze 
properties of $M_1(k)$ and $M_2(k)$.
For all $N\in \mathbb N$, we have
\begin{equation}\label{per_ine2}
\begin{aligned}
\sum_{k=0}^N M_1(k)
=\sum_{k=0}^N\sum_{i\in\mathcal V}
&2\mu^1_i \big[\alpha_{2k} \big(f_i(v_i(2k))-f_i(v_i(2k+1))\big)\\
&+(\alpha_{2k}-\alpha_{2k+1})
\big(f_i(v_i(2k+1))-f_i(x^*)\big)
\big].
\end{aligned}
\end{equation}
By (\ref{lip_ine}), we obtain
\begin{equation*}
\begin{aligned}
\alpha_{2k} &\big|f_i(v_i(2k)) -f_i(v_i(2k+1))\big|
\le \alpha_{2k}L|v_i(2k)-v_i(2k+1)| \\
&\le \alpha_{2k}L \big(|v_i(2k)-x_i(2k)| +|v_i(2k+1)-x_i(2k+1)|
+|x_i(2k)-x_i(2k+1)| \big).
\end{aligned}
\end{equation*}

It follows from Lemma \ref{cons:lem} (ii) that
$$\lim_{N\to \infty} \sum_{k=0}^N\sum_{i\in\mathcal V}
\alpha_{2k}L \big(|v_i(2k)-x_i(2k)| +|v_i(2k+1)-x_i(2k+1)|\big) < \infty.$$
According to (\ref{alg}) and (\ref{pro_ine1}), we have $$|x_i(2k)-x_i(2k+1)|\le |x_i(2k)-v_i(2k)|+|x_i(2k+1)-v_i(2k)|\le |x_i(2k)-v_i(2k)|+\alpha_{2k}L.$$
In view of Assumption \ref{step_ass}, 
$\lim_{N\to \infty} \sum_{k=0}^N \sum_{i\in\mathcal V} \alpha_{2k}|x_i(2k)-x_i(2k+1)|<\infty$.
As a result,
$$\lim_{N\to \infty} \sum_{k=0}^N\sum_{i\in\mathcal V} \mu^1_i \alpha_{2k} \big|f_i(v_i(2k))-f_i(v_i(2k+1))\big|<\infty.$$

Revisit the last term of (\ref{per_ine2}).
Notice that 
(\ref{alg}) implies that both $x_i(k)$ and $y(k)$ are bounded because $X$ is a compact set.
Due to the continuity of $f_i$, 
there must be a scalar $M_0$ such that $|f_i(v_i(k))-f_i(x^*)|\le M_0$.
Because $\alpha_k$ is non-increasing, we have
\begin{equation*}
\begin{aligned}
\sum_{k=0}^N\sum_{i\in\mathcal V}&
\mu^1_i (\alpha_{2k}-\alpha_{2k+1})
\big|f_i(v_i(2k+1))-f_i(x^*)\big| \le \sum_{k=0}^N(\alpha_{2k}-\alpha_{2k+1})M_0 \\
&\le (\alpha_0-\alpha_1)M_0 +\sum_{k=1}^N (\alpha_{2k-1}-\alpha_{2k+1})M_0
< \alpha_0 M_0 <\infty.
\end{aligned}
\end{equation*}

In summary, we conclude that $\lim_{N\to \infty}\sum_{k=0}^N M_1(k)$ is bounded.
By a similar procedure for discussing $M_1(k)$, we also have $\lim_{N\to \infty}\sum_{k=0}^N M_2(k) <\infty$.

\emph{Step 3.}
In this step, we show that $x_i(k)$ converges to a point in $X_p^*$.
Note that 
$$f_i(v_i(k))-f_i(x^*)=f_i(v_i(k))-f_i(y(k))
+f_i(y(k))-f_i(x^*).$$
By re-arranging the terms of (\ref{per_ine1}) and summing these relations over the time interval $k=0$ to $N$, we have
\begin{equation}\label{per_ine4}
\begin{aligned}
\chi_{N+1}
\le& \chi_0
+\sum_{k=0}^N (\alpha_{2k}^2+2\alpha_{2k+1}^2
+\alpha_{2k+2}^2) L^2
+\sum_{k=0}^N \big(M_1(k)+M_2(k)\big)\\
&\quad-\sum_{k=0}^{2N+1}\sum_{i\in\mathcal V}2 (\mu^1_i+\mu^2_i)\alpha_{k}
\big(f_i(v_i(k))-f_i(y(k))\big)\\
&\quad-\sum_{k=0}^{2N+1}\sum_{i\in\mathcal V}
2 (\mu^1_i+\mu^2_i)\alpha_{k}
\big(f_i(y(k))-f_i(x^*)\big).\\
\end{aligned}
\end{equation}
In light of (\ref{lip_ine}) and Lemma \ref{cons:lem} (ii), we obtain
\begin{equation*}
\begin{aligned}
\lim_{N\to\infty}
-\sum_{k=0}^{2N+1}&\sum_{i\in\mathcal V}
2 (\mu^1_i+\mu^2_i)\alpha_{k}
\big(f_i(v_i(k))-f_i(y(k))\big) \\
\le &\lim_{N\to\infty} \sum_{k=0}^{2N+1}\sum_{i\in\mathcal V}
2 (\mu^1_i+\mu^2_i)\alpha_{k}|v_i(k)-y(k)|
<\infty.
\end{aligned}
\end{equation*}

Due to the convexity and continuity of $f_i$, for any $x^*\in X_p^*$ and $\epsilon>0$ such that ${\rm dist}(y(k),X_p^*)>\epsilon$,
there exists $\delta>0$ such that $\sum_{i\in\mathcal V}
2(\mu_{1,i}+\mu_{2,i})\alpha_k
\big(f_i(y(k))-f_i(x^*)\big)>\delta$. 
Then 
$$\sum_{k=0}^{2N+1}\sum_{i\in\mathcal V}
2(\mu_{1,i}+\mu_{2,i})\alpha_k
\big(f_i(y(k))-f_i(x^*)\big)
\ge \sum_{k=0}^{2N+1}\sum_{i\in\mathcal V}
2\delta(\mu_{1,i}+\mu_{2,i})\alpha_k.$$
Under Assumption \ref{step_ass}, 
if ${\rm dist}(y(k),X_p^*)>\epsilon$ for any $\epsilon>0$,
the right hand of (\ref{per_ine4}) tends to $-\infty$ as $N$ tends to infinity. This contradicts with $\chi_{N+1}\ge 0$.
Therefore, we conclude that
$\lim_{k\to \infty}\inf {\rm dist}(y(k), X_p^*)=0$. Since $y(k)$ is bounded, there must be at least a limit point. In view of 
$\lim_{k\to\infty}\inf {\rm dist} (y(k),X_p^*)=0$, 
one of the limit points, denoted by $y^*$, must be in $X_p^*$.

It follows from (\ref{per_ine1}) that
\begin{equation*}
\begin{aligned}
\lim_{k\to\infty}\sup \chi_{k} 
& \le \lim_{k\to\infty} \inf\chi_{k}.
\end{aligned}
\end{equation*}
As a result, the scalar sequence $\chi_{k}$ 
is convergent. Then the limit point $y^*$ is unique.
Due to $\lim_{k\to \infty}|x_i(k)-y(k)|=0$, we conclude that the sequence $\{x_i(k)\}$ converges to the same point $y^*$ in $X_p^*$ .
This completes the proof.

\subsection{Proof of Theorem \ref{con:nonperiodic}}
\label{pf:noper}
Firstly, we show that  if $p=2$, (\ref{alg}) is convergent under quasi-periodic
graph sequences.
Here, notations are the same as that in the proof of Theorem \ref{con:periodic}.
Consider  $\mathcal G(k)$ switching between $\mathcal G_{A_1}$ and $\mathcal G_{A_2}$.
For quasi-periodic graphs, there are two cases:
$\mathcal G(2k)=\mathcal G_{A_1}$ and $\mathcal G(2k+1)=\mathcal G_{A_2}$;
$\mathcal G(2k)=\mathcal G_{A_2}$ and $\mathcal G(2k+1)=\mathcal G_{A_1}$.
If $\mathcal G(2k)=\mathcal G_{A_1}$ and $\mathcal G(2k+1)=\mathcal G_{A_2}$. As proved in (\ref{per_ine1}), we have
\begin{equation}\label{qua_ine1}
\begin{aligned}
\chi_{k+1} &\le \chi_k
+(\alpha_{2k}^2+2\alpha_{2k+1}^2
+\alpha_{2k+2}^2) L^2
-\sum_{i\in\mathcal V}2 (\mu^1_i+\mu^2_i)\alpha_{2k}
\big(f_i(v_i(2k))-f_i(x^*)\big) \\
&\quad-\sum_{i\in\mathcal V} 
2(\mu^1_i+\mu^2_i)\alpha_{2k+1}
\big(f_i(v_i(2k+1))-f_i(x^*)\big)
+M_1(k)+M_2(k).
\end{aligned}
\end{equation}
If $\mathcal G(2k)=\mathcal G_{A_2}$ and $\mathcal G(2k+1)=\mathcal G_{A_1}$, by a similar way for discussing $\chi_k$, we obtain
\begin{equation}\label{qua_ine2}
\begin{aligned}
\tilde \chi_{k+1} &\le \tilde \chi_k
+(\alpha_{2k}^2+2\alpha_{2k+1}^2
+\alpha_{2k+2}^2) L^2
-\sum_{i\in\mathcal V}2 (\mu^1_i+\mu^2_i)\alpha_{2k}
\big(f_i(v_i(2k))-f_i(x^*)\big) \\
&\quad-\sum_{i\in\mathcal V} 
2(\mu^1_i+\mu^2_i)\alpha_{2k+1}
\big(f_i(v_i(2k+1))-f_i(x^*)\big)
+M_3(k)+M_4(k),
\end{aligned}
\end{equation}
where $\tilde\chi_k=\sum_{i\in\mathcal V}
\mu^2_i|x_i(2k)-x^*|^2 
+\sum_{i\in\mathcal V} \mu^1_i|x_i(2k+1)-x^*|^2$,
$$M_3(k)=\sum_{i\in\mathcal V}2  
\mu^2_i\big[\alpha_{2k}
\big(f_i(v_i(2k))-f_i(x^*)\big)
-\alpha_{2k+1}
\big(f_i(v_i(2k+1))-f_i(x^*)\big)\big],$$
and moreover,
$$M_4(k)=\sum_{i\in\mathcal V} 
2\mu^1_i\big[ \alpha_{2k+1}
\big(f_i(v_i(2k+1))-f_i(x^*)\big)
-\alpha_{2k+2}
\big(f_i(v_i(2k+2))-f_i(x^*)\big)\big].$$

Define 
$\iota_k=\max\{\chi_k, \tilde \chi_k\}$,
$M_5(k)=\max\{M_1(k),M_3(k)\}$
and $M_6(k)=\\
\max\{M_2(k),M_4(k)\}$.
Note that $M_3(k)$ and $M_4(k)$ have similar properties as $M_1(k)$ proved in Section \ref{pf:con_per}.
Combining (\ref{qua_ine1}) and (\ref{qua_ine2}), we derive
\begin{equation}\label{qua_ine}
\begin{aligned}
\iota_{k+1} &\le \iota_k
+(\alpha_{2k}^2+2\alpha_{2k+1}^2
+\alpha_{2k+2}^2) L^2
-\sum_{i\in\mathcal V}2 (\mu^1_i+\mu^2_i)\alpha_{2k}
\big(f_i(v_i(2k))-f_i(x^*)\big) \\
&\quad-\sum_{i\in\mathcal V} 
2(\mu^1_i+\mu^2_i)\alpha_{2k+1}
\big(f_i(v_i(2k+1))-f_i(x^*)\big)
+M_5(k)+M_6(k).
\end{aligned}
\end{equation}

Clearly, 
$\lim_{N\to \infty}\sum_{k=0}^{\infty} \big(M_5(k)+M_6(k)\big)<\infty$.
By a similar procedure for discussing $\chi_k$ and $x_i(k)$ in  Section \ref{pf:con_per}, we can prove that $x_i(k)$ converges to a point in $X_p^*$ for all $i \in \mathcal V$.

In the following, 
we show that there exists a graph sequence such that (\ref{alg}) is not convergent if $p\ge3$. We consider $p=3$ for simplification, and note that the result can be easily extended to cases of $p>3$.
In view of Theorem \ref{con:periodic}, for the periodic graph sequence with the order $\mathcal G_{A_1}\rightarrow \mathcal G_{A_2}\rightarrow \mathcal G_{A_3}$, (\ref{alg}) converges to a point in $X^*_{sp}$, where $X^*_{sp}$ is the solution set of
\begin{equation*}\label{non_per1}
\min~~\sum_{i\in\mathcal V}
(\mu^1_i+\mu^2_i+\mu^3_i) f_i(x),~~~~
{\rm s.t.}~~x\in X,
\end{equation*}
where $\mu^1$, $\mu^2$ and $\mu^3$ are the Perron vectors of $A_3A_2A_1$, $A_2A_1A_3$ and $A_1A_3A_2$, respectively.
Similarly, for the periodic graph sequence with the order 
$\mathcal G_{A_1} \rightarrow\mathcal G_{A_3}\rightarrow \mathcal G_{A_2}$, 
(\ref{alg}) converges to a point in $\tilde X^*_{sp}$, where $\tilde X^*_{sp}$ is the solution set of
\begin{equation*}\label{non_per2}
\min~~\sum_{i\in\mathcal V}
(\tilde\mu^1_i+\tilde\mu^2_i+
\tilde \mu^3_i) f_i(x),~~~~
{\rm s.t.}~~x\in X,
\end{equation*}
where $\tilde\mu^1$, $\tilde\mu^2$ and $\tilde\mu^3$ are the Perron vectors of $A_2A_3A_1$, $A_3A_1A_2$, and $A_1A_2A_3$, respectively.
If $\bigcap_{i\in \mathcal V} X_i^*=\emptyset$, it follows from Remark \ref{con:suf_nec} that there exist $(\mu^1_i+\mu^2_i+\mu^3_i)$ and 
$(\tilde\mu^1_i+\tilde\mu^2_i+
\tilde \mu^3_i)$
such that $X^*_{sp}\bigcap \tilde X^*_{sp}=\emptyset$.
Because the weight matrices can be chosen freely under Assumptions \ref{gra_ass} and \ref{wei_ass}, there will always be 
$A_1$, $A_2$ and $A_3$ 
such that $X^*_{sp}\bigcap \tilde X^*_{sp}=\emptyset$.

By a similar way as the proof of Theorem \ref{non_con_the}, we can construct  time sequences $\{t_k\}$ and $\{s_k\}$,
and a graph sequence $\{\mathcal G(k)\}$ switching between $\mathcal G_{A_1}\rightarrow \mathcal G_{A_2}\rightarrow \mathcal G_{A_3}$ and $\mathcal G_{A_1} \rightarrow\mathcal G_{A_3}\rightarrow \mathcal G_{A_2}$ at time intervals $[3t,3(t+1))$. Then  (\ref{alg}) is not convergent under 
$\{\mathcal G(k)\}$. 
This completes the proof. 
$\hfill\Box$

\section{Numerical simulations}
Here, an illustrative example is carried out to verify the theoretical results presented in Section \ref{results}.

Similar to \cite{liu2017convergence}, we  employ (\ref{alg}) to solve the constrained LASSO (the least
absolute shrinkage and selection operator) regression problem, which is formulated as
\begin{equation*}\label{ex1}
\min~~\sum_{i\in \mathcal V} 
\big(\frac 12|x-q_i|^2+\sigma |x|_1 \big),
\quad {\rm s.t.}~~x'x \le 1,
\end{equation*}
where $q_i \in \mathbb R^4$ are known vectors,  $\sigma=0.1$ is a regularization parameter, and $x \in \mathbb R^4$ is the decision variable.
There are twenty agents in the multi-agent network, where agent $i$ only knows $q_i$.
Take $\alpha_k=k^{-0.6}$ and each entry of $x_i(0)$ from $[0,0.1]$ for (\ref{alg}).

Firstly, we show that
there generally exists a graph sequence such that (\ref{alg}) is not convergent.
We generate each entry of $q_i$ by a uniform distribution over $[-2,2]$, and two strongly connected graphs $\mathcal G_{A_1}$ and  $\mathcal G_{A_2}$ under Assumption \ref{wei_ass}.
Then we compute $X_{\mu(A_l)}^*$ under $\mathcal G(k) =\mathcal G_{A_l}$ for $l=\{1,2\}$. Following that, we construct a time sequence $\{t_k\}$, and a graph sequence $\mathcal G(k)$ by (\ref{graph_seq_cons}).

\begin{figure}[H]
\label{con_fig}
\centering
\includegraphics[scale=0.6]{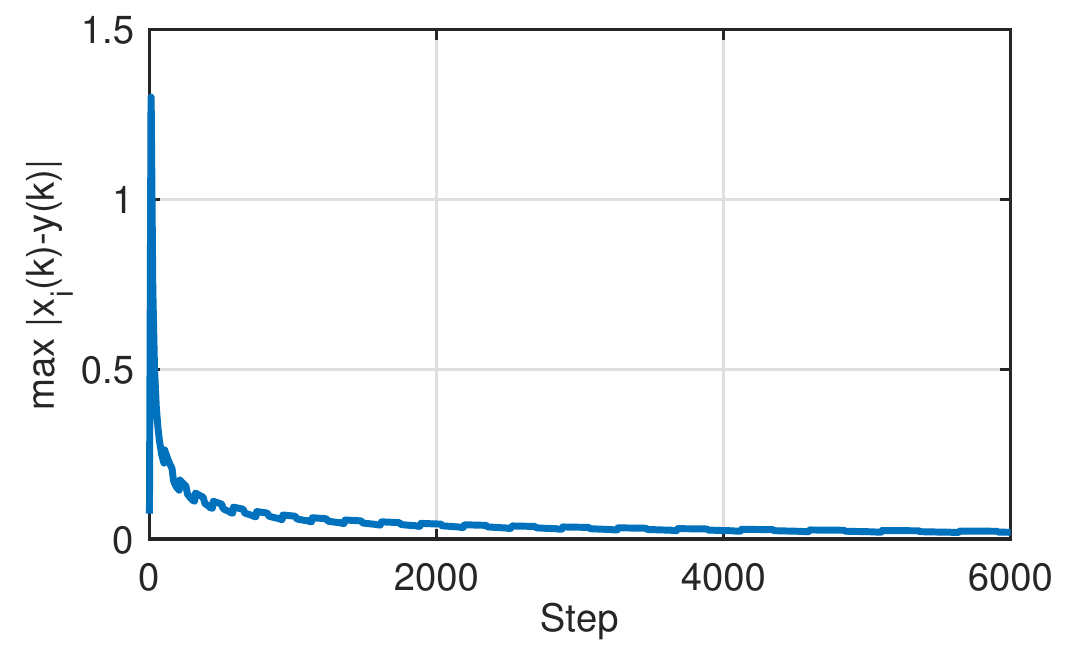}
\caption{The consensus result of (\ref{alg}).}
\end{figure}

\begin{figure}[H]
\label{non_con_fig}
\centering
\includegraphics[scale=0.6]{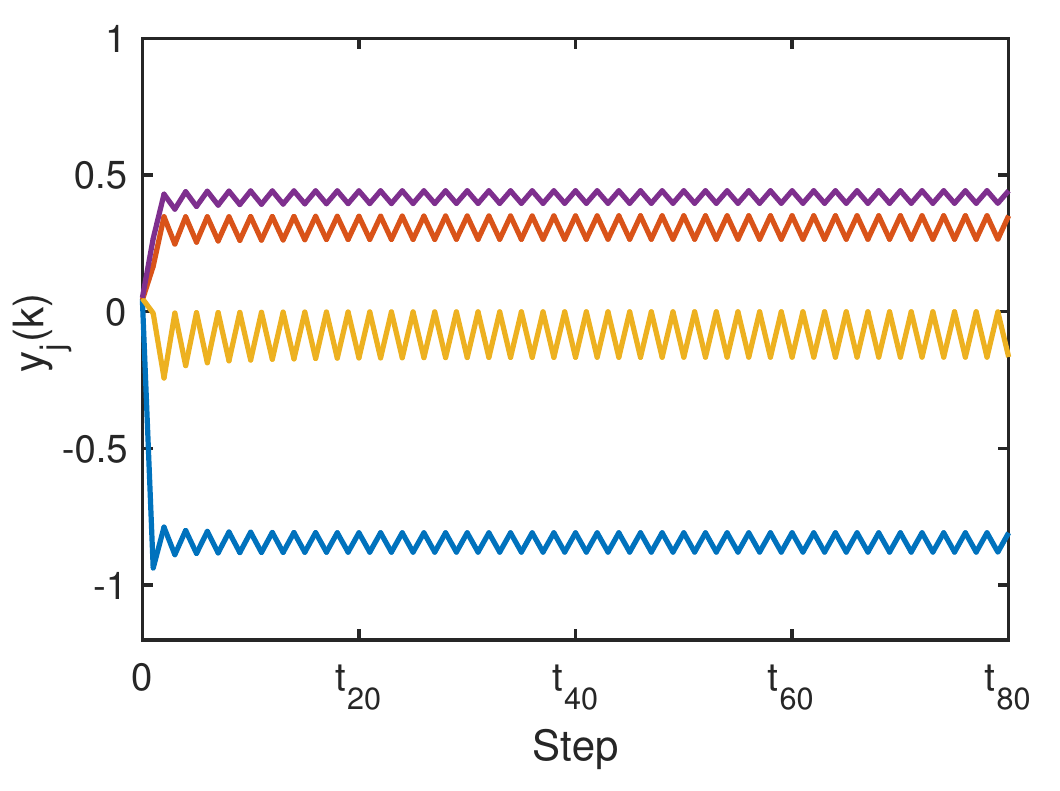}
\caption{The trajectories of $y_i(k)$ by  (\ref{alg}).}
\end{figure}

Figs. \ref{con_fig} and \ref{non_con_fig} show trajectories of $\max_{i\in\mathcal V}|x_i(k)-y(k)|$ and $y_j(k)$, respectively, where $y(k)=\frac 1n\sum_{i\in \mathcal V}  x_i(k)$ and $y_j(k)$ is the $j$-th entry of $y(k)$. It can be concluded that all agents achieve consensus because of  $\lim_{k\to \infty}|x_i(k)-y(k)|=0$. Furthermore,  $x_i(k)$ is not convergent due to the oscillation of $y(k)$.

\begin{figure}[H]
\label{con:suf_nec}
\centering
\includegraphics[scale=0.6]{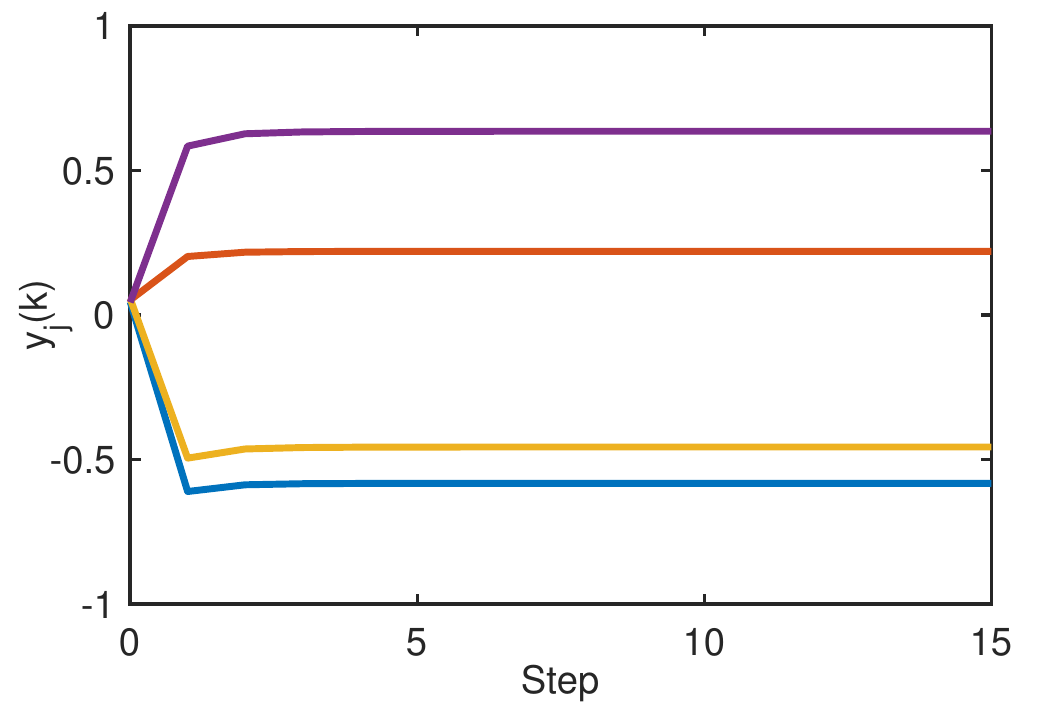}
\caption{The trajectories of $y_i(k)$ with $\bigcap_{i\in \mathcal V}X_i^* \not= \emptyset$.}
\end{figure}

Theorem \ref{con:nec_suf} is verified as follows. 
If the intersection of optimal solution sets to each agent is empty, then there is a graph sequence such that (\ref{alg}) is not convergent.
The result is very similar to Fig. \ref{non_con_fig}, and is omitted here.
Consider all agents being with the same $q_i$, and then $\bigcap_{i\in \mathcal V}X_i^* \not=\emptyset$. 
Fig. \ref{con:suf_nec} shows the trajectories of $y_i(k)$, where the network $\mathcal G(k)$ switches freely within four strongly connected graphs.
The result indicates the convergence of (\ref{alg}) in this case.

\begin{figure}[H]
\label{per:sim}
\centering
\includegraphics[scale=0.6]{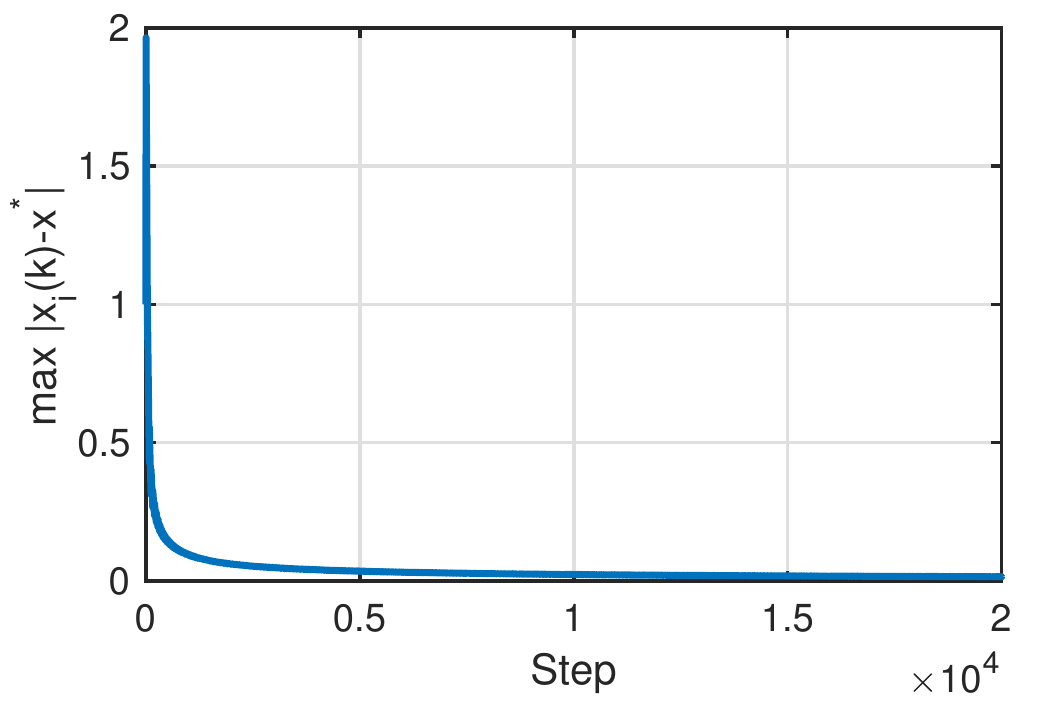}
\caption{The convergence result of algorithm (\ref{alg}) under a periodic graph sequence.}
\end{figure}

In the following, we demonstrate convergence results of (\ref{alg}) under periodic and quasi-periodic graph sequences. 
We generate each entry of $q_i$ by a uniform distribution over $[-2,2]$, and three graphs $\mathcal G_{A_1}$, $\mathcal G_{A_2}$ and $\mathcal G_{A_3}$, where $\mathcal G_{A_1} \cup \mathcal G_{A_2}$ is strongly connected.
Consider the graph sequence given by
$$\mathcal G_{A_1}\rightarrow \mathcal G_{A_2} \rightarrow \mathcal G_{A_3} \rightarrow \mathcal G_{A_1}\rightarrow \mathcal G_{A_2} \rightarrow \mathcal G_{A_3}\rightarrow \cdots.$$
By the centralized projected gradient algorithm, we compute the optimal solution $x^*$ to (\ref{per_for}).
Then we employ (\ref{alg}) for this problem, and plot the trajectory of $\max_{i\in \mathcal V}|x_i(k)-x^*|$ in Fig. \ref{per:sim}.
In this case, (\ref{alg}) achieves a solution to (\ref{per_for}) due to $\lim_{k\to\infty}|x_i(k)-x^*|=0$.

\begin{figure}[H]
\label{per_not_con}
\centering
\includegraphics[scale=0.6]{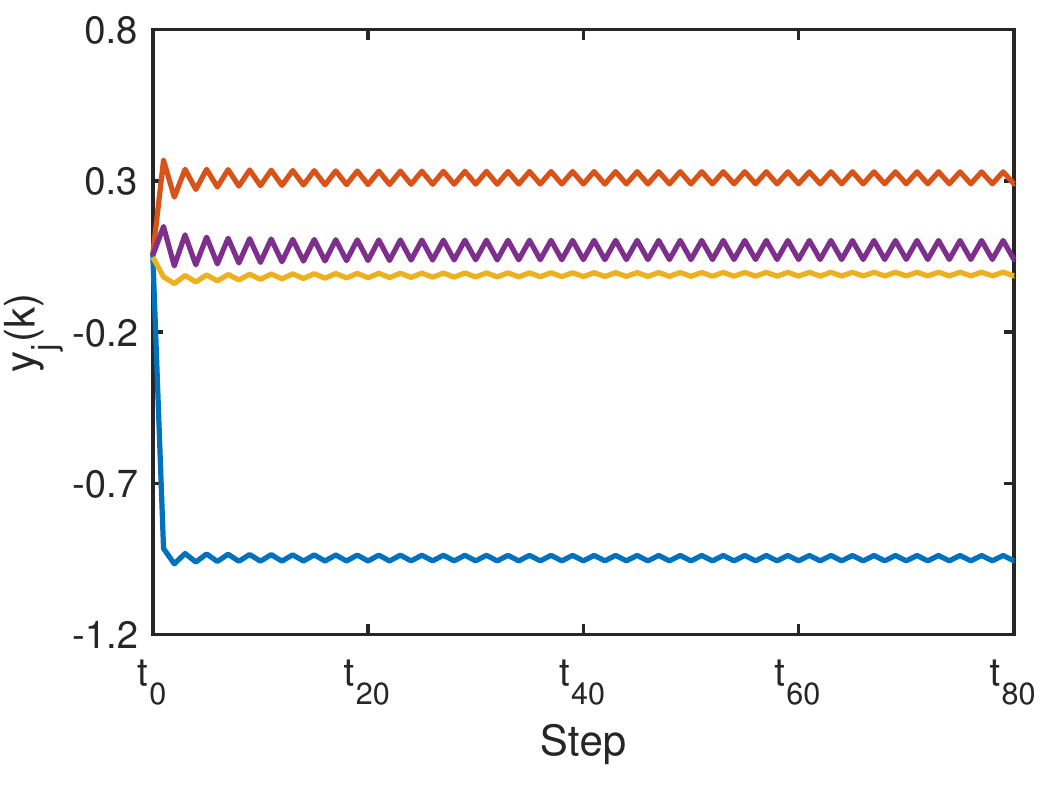}
\caption{The convergence result of (\ref{alg}) under a quasi-periodic graph sequence.}
\end{figure}

We discuss the case of quasi-periodic graph sequences. To be specific, $\mathcal G(k)$ switches between $\mathcal G_{A_1}\rightarrow \mathcal G_{A_2} \rightarrow \mathcal G_{A_3}$ and $\mathcal G_{A_3}\rightarrow \mathcal G_{A_2} \rightarrow \mathcal G_{A_1}$ at each time interval $[3t, 3(t+1))$.
By (\ref{graph_seq_cons}), we construct a graph sequence such that (\ref{alg}) is not convergent, and show the result in Fig. \ref{per_not_con}.

\begin{figure}[H]
\label{per_like_con}
\centering
\includegraphics[scale=0.6]{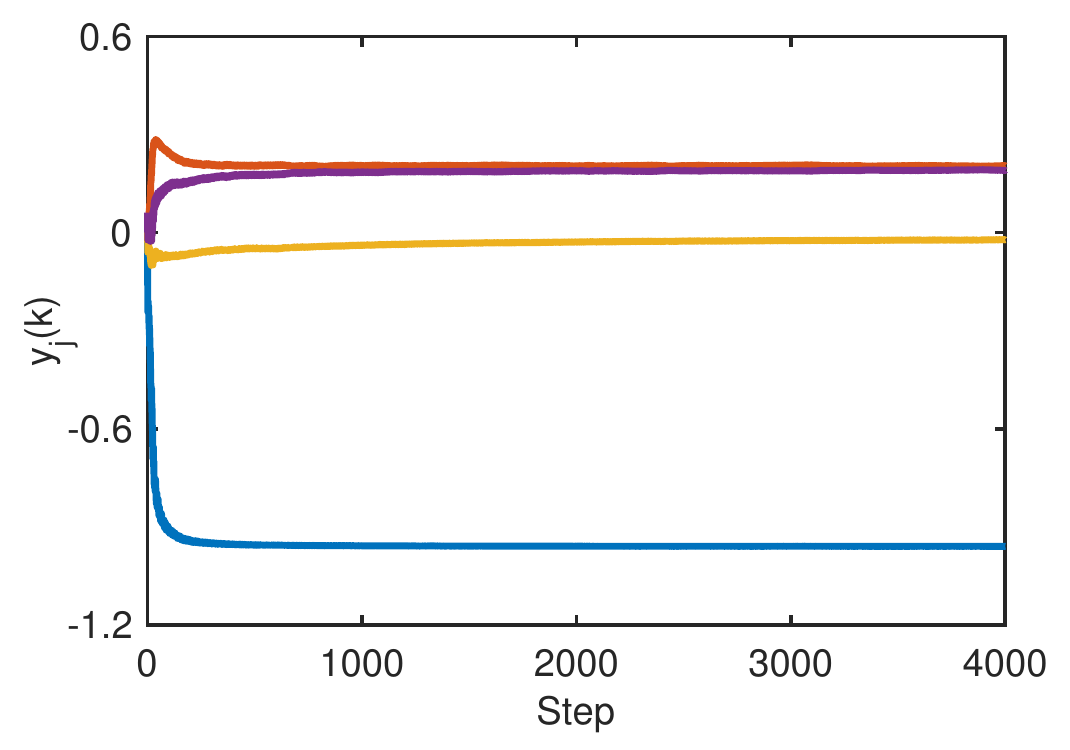}
\caption{The convergence result of (\ref{alg}) under a quasi-periodic graph.}
\end{figure}

Finally, we consider $\mathcal G(k)$ switching freely between $\mathcal G_{A_1}\rightarrow \mathcal G_{A_2}$ and $\mathcal G_{A_2} \rightarrow \mathcal G_{A_1}$ at each time interval $[2t, 2(t+1))$. Fig. \ref{per_like_con} shows the state trajectory of (\ref{alg}). The result indicates the convergence of (\ref{alg}) in this case.
Figs. \ref{per_not_con} and \ref{per_like_con} imply the correctness of Theorem \ref{con:nonperiodic}.

\section{Conclusions}

This paper aimed at investigating convergence properties of a distributed projected subgradient algorithm, where weight matrices of the time-varying communication network were only required to be row stochastic, i.e, the network might be weight-unbalanced. Firstly, it was proved that there generally existed a graph sequence such that the algorithm was not convergent if the network switched freely within finitely many general graphs.
Then  to guarantee the convergence of this algorithm for any uniformly strongly connected general graph sequence, it was provided a necessary and sufficient condition,  i.e., the intersection of optimal solution sets to all local optimization problems was not empty.
Following that, it was found that the algorithm was convergent under periodically switching general graph sequences, and optimized a weighted sum of local cost functions.
Furthermore, the periodic condition was slightly relaxed by quasi-periodic graph sequences, and it was shown that the algorithm was always convergent for any quasi-periodic graph sequence if and only if the network switched between two graphs.
Finally, numerical simulations were carried out for illustration.

\bibliographystyle{IEEEtran}
\bibliography{reference}

\end{document}